\documentclass[a4paper, twoside, 11pt, english]{article}
 
\usepackage{amsmath, amsthm, amsfonts, amssymb, mathabx}
\usepackage[T1]{fontenc}
\usepackage[utf8]{inputenc}
\usepackage{newtxtext, courier, euler}
\usepackage[scaled=0.95]{helvet}
\usepackage[cal=euler, scr=boondoxo, scrscaled=1.05]{mathalfa}
\usepackage{fancyhdr, titlesec, url, enumerate, setspace}
\usepackage{tikz}
\usepackage[all,knot,poly]{xy}
\usepackage{algorithm2e}
\usepackage{babel}
\usepackage{hyperref}

\titleformat{\section}[block]{\scshape\filcenter\LARGE}{\thesection.}{.5em}{}
\titleformat{\subsection}[runin]{\bfseries}{\thesubsection.}{.5em}{}[.]
\titlespacing{\subsection}{0pt}{\topsep}{.5em}

\newtheoremstyle{ntheorem}%
	{\topsep}{\topsep}{\itshape}{0pt}{\bfseries}{.}{.5em}%
	{\thmname{#1\hspace{.4em}}\thmnumber{#2}\thmnote{ (#3)}}
	
\newtheoremstyle{ndefinition}%
	{\topsep}{\topsep}{\normalfont}{0pt}{\bfseries}{.}{.5em}%
	{\thmname{#1\hspace{.4em}}\thmnumber{#2}.\thmnote{ (#3)}}
	
\newtheoremstyle{nremark}%
	{\topsep}{\topsep}{\normalfont}{0pt}{\itshape}{.}{.5em}%
	{\thmname{#1\hspace{.4em}}\thmnumber{#2}.\thmnote{ (#3)}}

\theoremstyle{ntheorem}
  	\newtheorem{theorem}[]{Theorem}

\theoremstyle{ndefinition}

\fancypagestyle{plain}{
\fancyhf{}\fancyfoot[LC,RC]{}
\fancyfoot[LE,RO]{$\thepage$}

}

\newcommand{\K}{\mathbb{K}}
\newcommand{\Nb}{\mathbb{N}}

\newcommand{\Mr}{\mathcal{M}}

\newcommand{\Ur}{\mathcal{U}}

\newcount\hh
\newcount\mm
\mm=\time
\hh=\time
\divide\hh by 60
\divide\mm by 60
\multiply\mm by 60
\mm=-\mm
\advance\mm by \time
\def\hhmm{\number\hh:\ifnum\mm<10{}0\fi\number\mm}

\def\cone{\mathrm{cone}}
\renewcommand{\div}[1]{\mathcal{#1}}

\newcommand{\comp}[1]{#1^{\perp}}
\newcommand{\compp}[2]{#2^{\perp(#1)}}

\def\mult{\mathrm{Mult}}
\def\cmult{\mathrm{^{\complement}Mult}}
\def\cnonmult{\mathrm{^{\complement}NMult}}
\def\nonmult{\mathrm{NMult}}

\DeclareMathOperator{\lm}{lm}

\SetAlgorithmName{Procedure}{List of procedures}

\DeclareMathOperator{\integralcond}{{\bf IntCond}}

\renewcommand{\leq}{\leqslant}
\renewcommand{\geq}{\geqslant}

\def\wo{\preccurlyeq_{wo}}
\def\wostrict{\prec_{wo}}

\newcommand{\auteur}[3]{
\noindent
\begin{minipage}[t]{.45\textwidth}
\begin{flushright}
\textsc{#1} \\
{\footnotesize\textsf{#2}}
\end{flushright} 
\end{minipage}
\qquad
\begin{minipage}[t]{.45\textwidth}
#3
\end{minipage}
}

\begin{document}
\thispagestyle{empty}

\begin{center}

\begin{huge}
{\scshape Maurice Janet's algorithms on \\ systems of linear partial differential equations}

\end{huge}

\bigskip
\hrule height 1.5pt 
\bigskip

\begin{Large}
{\scshape Kenji Iohara
\qquad Philippe Malbos}
\end{Large}

\vspace{2cm}

\begin{small}\begin{minipage}{14cm}
\noindent\textbf{Abstract --}
This article describes the emergence of formal methods in theory of partial differential equations (PDE) in the French school of mathematics through Janet's work in the period 1913-1930.
In his thesis and in a series of articles published during this period, Janet introduced an original formal approach to deal with the solvability of the problem of initial conditions for finite linear PDE systems. His constructions implicitly used an interpretation of a monomial PDE system as a generating family of a multiplicative set of monomials. He introduced an algorithmic method on multiplicative sets to compute compatibility conditions, and to study the problem of the existence and the uniqueness of a solution to a linear PDE system with given initial conditions. The compatibility conditions are formulated using a refinement of the division operation on monomials defined with respect to a partition of the set of variables into multiplicative and non-multiplicative variables.
Janet was a pioneer in the development of these algorithmic methods, and the completion procedure that he introduced on polynomials was the first one in a long and rich series of works on completion methods which appeared independently throughout the 20th century in various algebraic contexts.

\medskip


\medskip

\smallskip\noindent\textbf{M.S.C. 2010 -- Primary:} 01-08, 01A60. \textbf{Secondary:} 13P10, 12H05, 35A25.
\end{minipage}\end{small}

\vspace{2cm}

\begin{small}\begin{minipage}{12cm}
\renewcommand{\contentsname}{}
\setcounter{tocdepth}{1}
\tableofcontents
\end{minipage}
\end{small}
\end{center}

\clearpage

\section{Introduction}

This article presents the emergence of formal methods in the theory of partial differential equations, PDE for short, in the French school of mathematics through the works of the French mathematician Maurice Janet  in the period from 1913 to 1930. Janet was a very singular mathematician, who had been able to bring out original algebraic and algorithmic methods for the analysis of linear PDE systems. This original contribution of Janet is certainly due to his open-mindedness, as is made clear by his scientific visits to Germany, during a very complex political context in Europe with the events around the First World War.
In particular, this relationship with the Göttingen school led him to appropriate Hilbert's constructive ideas from \cite{Hilbert1890} in the algebraic analysis of polynomial systems.
In the continuation of the works of Charles Riquier and \'{E}tienne Delassus, he defended a \emph{Doctorat és Sciences Mathématiques}
\cite{Janet20PhD} in 1920, where he introduced an original formal approach to deal with the solvability of the problem of initial conditions  for finite linear PDE systems.

In this article we briefly survey the historical background of the contribution of Janet and we present its precursory ideas on the algebraic formulation of completion methods for polynomial systems applied to the problem of analytic solvability of PDE systems. 
Certainly influenced by the work of David Hilbert, \cite{Hilbert1890}, its construction implicitly used an interpretation of a monomial PDE system as a generating family of a multiplicative set of monomials. He introduced an algorithmic method on multiplicative sets to compute compatibility conditions, and to study the problem on the existence and the uniqueness of a solution of a linear PDE system with an initial condition. The compatibility conditions are formulated using a refinement of the division operation on monomials defined with respect to a partition of the set of variables into multiplicative and non-multiplicative variables.
We will explain how Janet's constructions were formulated in terms of polynomial systems, but without the modern language of ideals introduced simultaneously by Emmy Noether  in 1921~\cite{Noether1921}.
Janet was a pioneer in the development of these algorithmic methods, and the completion procedure that he introduced on polynomials was the first in a long and rich series of works on completion methods which appeared independently throughout the 20th century in various algebraic contexts.
In this article, we do not present the theory developed by other pioneers on the formal approaches in the analysis of linear PDE systems, in particular the work of Joseph Miller Thomas, \cite{Thomas37}.

\subsection{Mathematical context in Europe after first World War}
In an early stage of his career, Janet developed his mathematical project in the context of the first World War, which caused a complicated political period in Europe. This war, which in particular involved France and Germany, had profoundly affected the European mathematical community. We refer the reader to  \cite{AubinGispertGoldstein14} and \cite{AubinGoldstein14} for an exposition of the impact of this war on the activities of the mathematical community in Paris. This wartime followed a very active period for mathematics in Paris, and destroyed the dynamism of the French mathematical school. 
Indeed, many mathematicians were mobilized and the communications between France and other countries became difficult, especially between France and Germany, its main enemy. We refer the reader to~\cite{Mazliak13} which presents an edition of private notes written by Janet in the autumn of 1912 during his visit to Göttingen. In these notes, Janet revealed his views on the very complex political situation in Europe during this period.

The wartime created a very special situation for scientific collaborations between France and Germany. Indeed, some scientists expressed suspicions about the work of the enemy country's scientists. In particular, Charles Émile Picard, whose family was very badly affected by the war due to the death of three of his five children, published a very critical text on German science in 1916, \cite{Picard1916}. He wrote~{\cite[P.36]{Picard1916}}.
\begin{quote}
\emph{
C'est une tendance de la science allemande de poser a priori des notions et des concepts, et d'en suivre indéfiniment les conséquences, sans se soucier de leur accord avec le réel, et même en prenant plaisir à s'éloigner du sens commun. Que de travaux sur les géométries les plus bizarres et les symbolismes les plus étranges pourraient être cités ; ce sont des exercices de logique formelle, où n'apparaît aucun souci de
distinguer ce qui pourra être utile au développement ultérieur de la science mathématique.
}
\end{quote}

\begin{quote}
\emph{It is a tendency of German science to introduce notions and concepts, and to follow their consequences indefinitely, without worrying about their agreement with reality, and even taking pleasure in departing  from common sense. How many works on the strangest geometries and the strangest symbolisms could be cited; they are exercises in formal logic, where there is no concern to distinguish what may be useful for the further developments of mathematical science?
}
\end{quote}
During this period, Picard had a significant influence on the French school of analysis. Consquently, such a strong position towards German scientists reflects the atmosphere of the period during which Janet was conducting his thesis work.
Nevertheless, having visited Germany, Janet had privileged relations with the German mathematical community. Janet's visit to Göttingen was described with details in~\cite{Mazliak13}. 
His work on formal methods for the solvability of linear PDE systems was influenced by the algebraic formalism developed during that period in Germany to deal with  finiteness problems in polynomial rings. Indeed, since Hilbert's seminal article, \cite{Hilbert1890}, these questions have been at the center of many works in Germany. It was in 1921 that the algebraic structure of ideals emerged, and the Noetherian property was clearly formulated. This was after a long series of works carried out by the German school, with the major contributions by D. Hilbert, \cite{Hilbert1890}, Richard Dedekind, and finally E. Noether.

In France, the formalist approach was not as well developed as in Germany and the reference text-books in algebra remained the great classics of the 19th century on the analysis of algebraic equation systems. In the 1920s, the main references were the book of Camille Jordan  on what he called substitution groups of algebraic equations, \cite{Jordan1870}, and the lectures on higher algebra by Joseph-Alfred Serret, \cite{Serret1849}. The book of J.-A. Serret had a great influence and was re-edited many times until 1928, \cite{Serret1928}.

\subsection{Maurice Léopold René Janet}
Janet was born on the 24th of October 1888 in Grenoble. He was raised in a family of six children that belonged to the French intellectual bourgeoisie. He entered the science section of the \emph{\'Ecole normale supérieure} in Paris in 1907. 
Jean-Gaston Darboux, Édouard Goursat, \'E. Picard and Jacques Hadamard  were among his teachers. In September 1912 he made a trip to Göttingen in Germany for a few months. This stay in Göttingen was thought to have been of great importance in the mathematical training of Janet. We refer the reader to \cite{Mazliak13} for more details on his travels to the University of Göttingen. 
He found a very rich intellectual community there and had many exchanges, both with foreign students visiting Göttingen like him: George Pólya, Lucien Godeaux, Marcel Riesz), and with prestigious teachers, Constantin Carathéodory, Richard Courant, Edmund  Landau, D. Hilbert and Felix Klein.\cite{Mazliak13}.
He also met Max Noether  and E. Noether, \cite{Mazliak13}.

The first two publications, \cite{Janet13a, Janet13b}, of Janet appeared in \emph{Comptes-rendus de l'Académie des sciences} in 1913 and deal with the analysis of PDE systems. The second publication \cite{Janet13b} concerns a generalization of Cauchy-Kowalevsky theorem under the formulation given by Ch. Riquier in \cite{Riquier10}.
While a lecturer at the University of Grenoble, he defended his thesis on the analysis of PDE systems, entitled  <<\emph{Sur les syst{\`e}mes d'{\'e}quations aux d{\'e}riv{\'e}es partielles}>>,  on June 26, 1920 at the Sorbonne in Paris, \cite{Janet20PhD}. The jury was composed of Gabriel Xavier Paul Koenigs, \'E. Goursat, \'Elie Cartan  and J. Hadamard.
In the preamble of his thesis dissertation, he payed a respectful tribute (\emph{Hommage respectueux et reconnaissant}) to \'Edouard Goursat and J. Hadamard.
Some results of his thesis were published in \emph{Journal de mathématiques pures et appliquées} in 1920, \cite{Janet20}.

Janet was promoted to Professor in 1920 in Nancy, then in Rennes in 1921. He became a  Professor in Caen in 1924 when Ch. Riquier retired from the University  and became Professor Emeritus. 
Finally, he became a Professor at the Sorbonne in Paris in 1945.
He was an invited speaker at the \emph{International Congress of Mathematicians} on three occasions: Toronto in 1924, \cite{Janet-ICM24}, Z\"{u}rich in 1932, \cite{Janet-ICM32} and Oslo in 1935, \cite{Janet-ICM36}. He was President of the \emph{Société Mathématiques de France} in 1948. He died in 1983.

\subsection{Formal methods in commutative algebras throughout the 20th century}
Most of the formal computational methods in commutative algebra and algebraic geometry developed throughout the 20th century were founded on extensive works in elimination theory in the period 1880-1915. As early as 1882, Kronecker \cite{Kronecker82} introduced multivariate resultants providing complete elimination methods for systems of polynomial equations.  Elimination theory culminated with the works of Julius K{\"o}nig (1849-1913) \cite{Koenig1903}, and of Macaulay \cite{Macaulay1903,Macaulay13,Macaulay16}. For an overview of the works of this period, the reader may consult an important book on algebra \cite{van-der-Waerden1930} written by Bartel Leendert van der Waerden based on lectures by E. Noether and E. Artin.

Independently, a computational approach to elimination in commutative algebra that consists to define a polynomial ideal throughout a generating family satisfying nice computational properties stated with respect to a monomial order appeared in different forms and in various contexts from the early 20th century. The first algebraic constructions using such a computational method appeared in \cite{Dickson13, Gunther41, Macaulay16}. 
Fifty years later, the notion of generating set of an ideal satisfying computational properties with respect to a monomial order appeared in the terminology of \emph{standard bases} in \cite{Hironaka64} for power series rings by Hironaka. In the same period, Bruno Buchberger (1942-) developed algorithmic approaches for commutative polynomial algebras, with effective constructions and a completion algorithm for calculating \emph{Gröbner bases}, \cite{Buchberger65}. Similar approaches were developed for non commutative algebras in \cite{Shirshov62, Bergman78}. 
Thereafter, developments of the theory of Gröbner bases has mainly been motivated by algorithmic problems such as computations with ideals, manipulating algebraic equations, computing linear bases for algebras, Hilbert series, and homological invariants.

Forty years before the work of B. Buchberger, Janet introduced algorithmic approaches to the completion of a generating family of a polynomial ideal into a generating family satisfying computational properties quite similar to Gröbner bases. As we explain in the following sections, the completion methods constitute the essential part of the theory developed by Janet. He introduced a procedure to compute a family of generators of an ideal having the involutive property, and called \emph{involutive bases} in the modern language. This property is used to obtain a normal form of a linear partial differential equation system.

Janet's procedure of computation of involutive bases used a refinement of the classical polynomial division, called \emph{involutive division}, which is appropriate to the reduction of linear PDE systems. 
The completion procedure that he introduced is quite similar to the one defined with respect to classical division by B. Buchberger in \cite{Buchberger65} to produce Gröbner bases.
Subsequently, another approach to the reduction of linear PDE systems by involutive divisions was introduced by J. M. Thomas, \cite{Thomas37}. The terminology \emph{involutive} first appeared in \cite{Gerdt97}. We refer the reader to \cite{Mansfield96} for a discussion on relation between this notion and that of involutivity in the work of \'E. Cartan. We refer also to \cite{Seiler10} for a complete account of algebraic involutivity theory. Finally, note that the work of Janet was forgotten for about a half-century, and was rediscovered by F. Schwarz in 1992 in \cite{Schwarz92}.

\subsection{Conventions and notations}
In order to facilitate the reading of the different mathematical constructions extracted from the publications of Janet, we have chosen to use modern mathematical formulations. 
We provide here, a dictionary between the terminology used by Janet and the terminology used nowadays in the theory of partial differential equations.

\medskip

\begin{center}
  \begin{tabular}{|c|c|c|}
\hline
\begin{tabular}{c}
Janet terminology\\
(in French) 
\end{tabular}
& Modern terminology &  \begin{tabular}{c} Subsection \\ in the article
\end{tabular}\\ 
\hline\hline
\emph{module de monômes} &  multiplicative cone & \ref{SS:MultiplicativeCone}\\
\hline
\emph{forme} &  homogeneous polynomial & \ref{SS:OnNotionOfModule}\\
\hline
\emph{module de formes} &  polynomial ideal & \ref{SS:OnNotionOfModule}\\
\hline
\emph{famille de monômes de type fini}$^*$ & Noetherian property & \ref{SS:OnAlgebraicFinitenessProperties}\\
\hline
\emph{système de cotes} & monomial order& \ref{SS:NotionOfCote}, \ref{SS:ParametricPrincipalJanet}\\
\hline
\emph{postulation} & coefficients of Hilbert's series &\ref{SS:OnNotionOfModule}\\
\hline
  \end{tabular}
\end{center}

\bigskip

\noindent ($^*$) Janet did not use any specific terminology, but formulated the notion as follows, {\cite[pp. 11]{Janet29}}:

\noindent {\bf Théorème.} - \emph{Une suite de monomes $M_1$, $M_2$, $\ldots$ telle que chacun d'eux n'est multiple d'aucun des précédents ne comprend qu'un nombre fini de monomes}.

\medskip

\noindent {\bf Theorem.} - \emph{A sequence of monomials $M_1$, $M_2$, $\ldots$ such that each monomial is not a multiple of any preceding one contains only a finite number of monomials.}

\medskip

The following notations will be used in this article.
The set of non-negative integers is denoted by~$\mathbb{N}$.
The polynomial ring on the variables $x_1,\ldots,x_n$ over a field $\K$ of characteristic zero is denoted by $\K[x_1,\ldots,x_n]$. A polynomial is either zero or it can be written as a sum of a finite number of non-zero \emph{terms}, each term being the product of a scalar in $\K$ and a \emph{monomial}.
We will denote by $\Mr(x_1,\ldots,x_n)$ the set of monomials in the ring $\K[x_1,\ldots,x_n]$. For a subset $I$ of $\{x_1,\ldots,x_n\}$ we will denote by $\Mr(I)$ the set of monomials in $\Mr(x_1,\ldots,x_n)$ whose variables lie in $I$. A monomial $u$ in $\Mr(x_1,\ldots,x_n)$ is written as $u=x_1^{\alpha_1}\ldots x_n^{\alpha_n}$, were the $\alpha_i$ are non-negative integers. The integer $\alpha_i$ is called the \emph{degree} of the variable $x_i$ in $u$, it will be also denoted by $\deg_i(u)$. For $\alpha=(\alpha_1,\ldots,\alpha_n)$ in $\mathbb{N}^n$, we denote $x^\alpha=x_1^{\alpha_1}\ldots x_n^{\alpha_n}$ and $|\alpha|=\alpha_1+\ldots+\alpha_n$.

For a set $\Ur$ of monomials of $\Mr(x_1,\ldots,x_n)$ and $1\leq i\leq n$, we denote by $\deg_i(\Ur)$ the largest possible degree in the variable $x_i$ of the monomials in $\Ur$, that is
\[
\deg_i(\Ur) \: = \:
\max\big(\deg_i(u)\;|\;u\in \Ur\,\big).
\]
We call the \emph{cone} of the set $\Ur$ the set of all multiples of monomials in $\Ur$, defined by
\[
\cone(\Ur) \: = \: \bigcup_{u\in \Ur} u\Mr(x_1,\ldots,x_n) \: = \:
\{\, uv \;|\; u \in \Ur,\; v \in \Mr(x_1,\ldots, x_n) \,\}.
\]
Finally, to a monomial $x^\alpha = x^{\alpha_1}x^{\alpha_2}\ldots x_n^{\alpha_n}$ we will associate the differential operator:
\[
D^\alpha \: = \: 
\frac{\partial^{|\alpha|}\;\;}{\partial x_1^{\alpha_1}\partial x_2^{\alpha_2}\ldots \partial x_n^{\alpha_n}}.
\]

\subsection*{Acknowledgements} 
We wish to thank Cameron Calk for his careful reading and useful remarks to improve this article, and the anonymous referees for fruitful comments about this article.

\section{Historical context of Janet's work}
\label{S:HistoricalContext}

Janet's contribution discussed in this article is part of a long series of works on partial differential equation systems. In order to introduce the motivations of Janet's results, this section outlines the main contributions on the study of systems of partial differential equations achieved in the 19th century. We present the historical background of exterior differential systems and of the questions on PDE. For a deeper discussion of the theory of differential equations and the Pfaff problem, we refer the reader to~\cite{Forsyth90, Weber00} or~\cite{Cartan1899}.

\subsection{Pfaff's problem}
\label{SS:PfaffProblem}
Motivated by problems in analytical mechanics\footnote{We refer the reader to \cite{Dugas50} concerning history of mechanical problems.}, Leonhard Euler  and Joseph-Louis Lagrange  initiated the so-called \emph{variational calculus}, cf. \cite{Lagrange88}, which led to the problem of solving first-order PDE. This theory serves as a guide to the Janet contributions. 
In 1772, J.-L. Lagrange considered in \cite{Lagrange72} a PDE of the following form 
\begin{equation}
\label{Equation:Lagrange0}
F(x,y,\varphi,p,q)=0
\quad\text{with}\quad 
p=\frac{\partial \varphi}{\partial x}
\quad\text{and}\quad
q=\frac{\partial \varphi}{\partial y},
\end{equation}
\emph{i.e.} a PDE of one unknown function $\varphi$ in two variables $x$ and $y$.  Lagrange's method to solve this PDE can be summarized in three steps as follows:
\begin{enumerate}[{\bf i)}]
\item Express the PDE (\ref{Equation:Lagrange0}) in the form 
\begin{equation}
\label{Equation:Lagrange}
q=F_1(x,y,\varphi,p)
\quad\text{with}\quad
p=\frac{\partial \varphi}{\partial x}
\quad\text{and}\quad q=\frac{\partial \varphi}{\partial y}.
\end{equation}
\item Forgetting  the fact that $p=\frac{\partial \varphi}{\partial x}$, we consider the following $1$-form
\[ \Omega=d\varphi-pdx-qdy=d\varphi-pdx-F_1(x,y,\varphi,p)dy,
\]
by regarding $p$ as some (not yet fixed) function of $x,y$ and $\varphi$.
\item If there exist functions $M$ and $\Phi$ in variables $x,y$ and $\varphi$ satisfying $M\Omega=d\Phi$, then $\Phi(x,y,\varphi)=C$ for some constant $C$. Solving this new equation, we obtain a solution $\varphi=\psi(x,y,C)$ to the given equation~(\ref{Equation:Lagrange}).
\end{enumerate}
In 1814-15, Johann Friedrich Pfaff  \cite{Pfaff15} treated the case of a PDE of one unknown function in $n$ variables, depending on the case when $n$ is even or odd.

Recall that any PDE of any order is equivalent to a first order PDE system, that is involving only first partial derivatives of the unknown functions. Thus, we exclusively consider systems of first order PDE with $m$ unknown functions of the form
\[
F_k\big( x_1, \ldots, x_n, \varphi^1, \ldots, \varphi^m, \frac{\partial \varphi^a}{\partial x_i}~(1\leq a \leq m, 1\leq i\leq n)\big)=0,
\quad\text{for}\quad 
1\leq k \leq r.
\]
Introducing the new variables $p^{a}_i$, the system is defined on the space with coordinates $(x_i, \varphi^{a},p_i^{a})$ and is given by
\[ \begin{cases} F_k(x_i, \varphi^{a},p_i^{a})=0, & \\
                         d\varphi^{a}- \displaystyle{\sum_{i=1}^n} p_i^{a}dx_i=0, & \\
                         dx_1 \wedge \ldots \wedge dx_n \neq 0. & \end{cases}
\]
Notice that the last condition means that the variables $x_1,\ldots , x_n$ are independent. Such a system is called a \emph{Pfaffian system}. 
One is interested in whether this system admits a solution or not, and whether or not a solution is unique under some conditions. These questions are \emph{Pfaff's problems}. An approach using differential invariants was one of the key ideas developed, in particular, by Sophus Lie  (cf. \cite{Lie84}),  Darboux (cf. \cite{Darboux82}), and Georg Frobenius (cf. \cite{Frobenius77}) etc. before \'E.~Cartan \cite{Cartan1899}. 
See, e.g.,  \cite{Forsyth90} and \cite{Hawkins05} for historical foundations on Pfaff's problems.

\subsection{Cauchy-Kowalevsky's theorem}
\label{Subsection:CauchyKowalevskyTheorem}
A naive approach to Pfaff's problems, having applications to mechanics in mind, is the question of the initial conditions.  In series of articles published in 1842, Augustin Louis  Cauchy  studied systems of PDE of first order in the following form: 
\begin{equation}
\label{Equation:Cauchy}
\frac{\partial \varphi^{a}}{\partial t}\: = \:f_a( t, x_1,\cdots, x_n)
+\sum_{b=1}^{m}\sum_{i=1}^{n} f_{a,b}^{i}( t, x_1,\ldots, x_n) \frac{\partial \varphi^{b}}{\partial x_i},
\quad\text{for}\quad 1\leq a\leq m, 
\end{equation}
where $f_a, f_{a,b}^i$ and $\varphi^1, \ldots, \varphi^m$ are functions of the $n+1$ variables   $t,x_1,\ldots,x_n$.
He showed that under the hypothesis of analyticity of the coefficients, the PDE system (\ref{Equation:Cauchy}) admits a unique analytic local  solution satisfying a given initial condition.

Sophie Kowalevsky  in 1875 considered PDE systems of the form
\begin{equation}
\label{Equation:Kowalevsky}
\frac{\partial^{r_a} \varphi^{a}}{\partial t^{r_a}}\: = \:\sum_{b=1}^{m}
\sum_{\substack{j=0 \\ j+\vert \alpha\vert \leq r_a}}^{r_a-1} f_{a,b}^{j,\alpha}( t, x_1,\ldots, x_n) \frac{\partial^{j+\vert \alpha\vert} \varphi^{b}}{\partial t^j \partial x^{\alpha}},
\end{equation}
for some $r_a \in \mathbb{Z}_{>0}$, $1\leq a\leq m$, and where $f_{a,b}^{j,\alpha}$ and $\varphi^1, \ldots, \varphi^m$ are functions of the $n+1$ variables $t,x_1,\ldots,x_n$, and $\alpha=(\alpha_1,\cdots, \alpha_n)$ in $(\mathbb{Z}_{\geq 0})^n$, with $\partial x^\alpha=\partial x_1^{\alpha_1} \ldots \partial x_n^{\alpha_n}$.
She proved in \cite{Kowalevsky75} that under the hypothesis of analyticity of the coefficients, the system (\ref{Equation:Kowalevsky}) admits a unique analytic local  solution satisfying a given initial condition. This result is now called the \emph{Cauchy-Kowalevsky theorem}. In her article \cite{Kowalevsky75}, she suspected that the form she has obtained was the normal form of any PDE system. However, she had no proof of this statement. She wrote, {\cite[pp. 24-25]{Kowalevsky75}}:
\begin{quote}
\emph{
Was dagegen die zweite Bedingung angeht, so bleibt allerdings noch zu untersuchen, ob ein Gleichungssystem von nicht normaler Form stets durch ein ähnliches Verfahren, wie es Jacobi bei einem System gewöhnlicher Differentialgleichungen angewandt hat, auf ein normales zurückgeführt werden könne, worauf ich aber hier nicht eingehen kann.
}
\end{quote}

\begin{quote}
\emph{
Regarding the second condition it remains to study the question of whether a system of equations not in normal form may always be reduced to a normal one by methods similar to the ones used by Jacobi for systems of ordinary differential equations, which I cannot go into here. 
}
\end{quote}

In his thesis in 1891, \cite{Bourlet91}, Charles Bourlet showed that any PDE system can be transformed into an equivalent PDE system of first order and proposed a notion of canonical form for such a system. He showed that for a completely integrable system, there is an analytic solution. He also showed that the normal form \eqref{Equation:Kowalevsky} due to S. Kowalevsky is not completely general by providing an example of a PDE system of one unknown function depending on the two independent variables.
Thus, finding a canonical form of more general PDE systems became an important problem in the context of generalizing the Cauchy-Kowalevsky theorem. In {\cite[\textsection 17]{Bourlet91}}, Bourlet wrote
\begin{quote}
\emph{
Ceci nous prouve que le théorème de M$\,^{\text{me}}$ de Kowalewski ne démontre pas l'existence des intégrales dans tous les cas où, dans le système à intégrer, le nombre des équations est égal au nombre des fonctions inconnues. Dans son Mémoire (Journal de Crelle, t. 80, pp. 25) M$\,^{\text{me}}$ de Kowalewski suppose que cette transformation soit possible en faisant, d'ailleurs, remarquer qu'elle ne peut assurer que cela soit toujours possible.
}
\end{quote}

\begin{quote}
\emph{ This shows that Mme. de Kowalewski's theorem does not establish the existence of the integrals in all the cases when, in the system to integrate, the number of equations and the number of unknown functions are equal. In her article (Crelle, vol. 80, pp. 25), Mme. de Kowalewski supposes that such a transformation is possible, while making a remark that, in any case, she cannot ensure that this is always possible.  
}
\end{quote}

The generalization of the Cauchy-Kowalevsky theorem to wider classes of linear PDE systems was at the origin of the works of C. Méray, \'E. Delassus, Ch. Riquier as explained in the next section. 
It was Janet who obtained a computational method to reach normal form of linear PDE systems for a class of systems satisfying a reducibility property recalled in Section~\ref{SS:NotionCanonicalPDESystems}.

\subsection{Grassman's differential rule} 
In 1844, Hermann Günther Grassmann  exhibited the rules of the exterior algebra computation in his book \cite{Grassmann44} on linear algebra, that is a relation of the type
\[
ab=-ba.
\]
Although this kind of relation was implicitly used in the computation of the determinant of a square matrix, as in a work of Carl  Jacobi  (cf. \cite{Jacobi27b} etc.), this approach was too abstract for the first half of the 19th century.

Cesare Burali-Forti  extensively applied this Grassmann's rule to elementary Geometry, in \cite{BuraliForti1897}, but had not treated what are now called differential forms. 
It was \'E. Cartan in 1899 \cite{Cartan1899} who introduced Grassmann's rule in differential calculus. This algebraic calculus allowed him to describe a PDE system by an exterior differential system which is independent of the choice of coordinates. This led to the so called \emph{Cartan-Kähler theory}, which is another motivation for the formal methods introduced by Janet for analysis on linear PDE systems. 
We refer the reader to \cite{Katz85} for the impact, in many fields of mathematics, of the introduction of Cartan's differential forms.
See, \cite{Kahler34} and \cite{Cogliati11} for historical accounts of the Cartan-K\"{a}hler theory and \cite{Griffiths83} and \cite{BC3G91} for exposition of this theory in modern language.

\section{Emergence of formal methods for linear PDE systems}
\label{S:EmergenceFormalMethodsPDE}

The Cauchy-Kowalevsky theorem gives conditions for the existence of solutions of the PDE system defined by~(\ref{Equation:Kowalevsky}) and satisfying some initial conditions.
Generalizations of this result to wider classes of linear PDE systems were investigated in France by Charles M\'{e}ray  Ch. Riquier 
and \'{E}. Delassus during the period 1880-1900. 
The first works in this direction seem to be those of a collaboration between Ch. Méray and Ch. Riquier,~\cite{MerayRiquier1889,MerayRiquier1890}. 
In the first of a series of three articles on the subject, \cite{Riquier93}, Ch. Riquier noted that a very small number of authors had, at that time, addressed the existence of integrals in a differential system involving any number of unknown functions and independent variables.

\subsection{Principal and parametric derivatives}
In the beginning of 1890s, following a collaboration with Ch.~M\'{e}ray, Ch.~Riquier initiated his research on finding normal forms of systems of (infinitely many) PDE for finitely many unknown functions with finitely many independent variables.
Ch.~Méray and Ch.~Riquier in \cite{MerayRiquier1890} analyzed S. Kowalevsky's proof in \cite{Kowalevsky75} with the objective of reducing a PDE to some notion of normal form.
It may be regarded as the first algorithmic method applied to the analysis of PDE systems.
They introduced the concept of \emph{principal} and \emph{parametric derivatives}, allowing them to make inductive arguments on sets of derivatives without having an explicit total order on these sets. They formulated this notion as follows,  {\cite[\textsection 2]{MerayRiquier1890}}:
\begin{quote}
\emph{
Dans un système d'équations différentielles partielles, il y a, relativement à chaque fonction inconnue, une distinction essentielle à faire entre les diverses variables indépendantes.
Nous appellerons variables principales d'une fonction inconnue déterminée celles par rapport auxquelles sont prises les dérivées de cette fonction qui constituent dans le Tableau du système les premiers membres des équations de la colonne correspondante. Pour la même fonction, toutes les autres variables seront paramétriques. 
}
\end{quote}

\begin{quote}
\emph{ In a system of partial differential equations, for each unknown function, there is, for each unknown function, an essential distinction to make between the various independent variables. We shall call the principal variables of an unknown function determined with respect to which the derivatives of this function are taken, that form the first members of the equations in the corresponding column in the table of the system. For the same function, all of the other variables are parametric. 
}
\end{quote}

The notions of principal and parametric derivatives as appearing in Méray-Riquier's work were not formally exposed in~\cite{MerayRiquier1890}. These notions would be formalized later by Janet in the elaboration of an algorithmic process for the computation of the normal form of a linear PDE system. We will present the Janet formulation of these derivatives in Section~\ref{S:JanetMonomialOrderOnDerivative}.

\subsection{The notion of cote}
\label{SS:NotionOfCote}
Ch.~Riquier noted in \cite{Riquier93} that the computation of normal forms for a PDE system requires defining a total order on the derivatives.
In this direction, he introduced the notion of \emph{cote} on derivatives in {\cite[pp. 66-67]{Riquier93}}, the first of a series of three articles published in a same volume of \emph{Annales Scientifiques de l'\'Ecole Normale Supérieure}.

\begin{quote}
\emph{
Désignant par
}

\medskip

\emph{
(1)\hspace{6cm} $x$, $y$, $\ldots$
}

\medskip

\emph{
les variables indépendantes, et par
}

\medskip

\emph{
(2)\hspace{6cm} $u$, $r$, $\ldots$
}

\medskip

\emph{
les fonctions inconnues d'un système différentiel quelconque, faisons correspondre à chacune des quantités (1), (2) $p$ entiers, positifs, nuls ou négatifs, que nous nommerons respectivement cote première, cote seconde, ..., cote $p^{ième}$ de cette quantité. Considérant ensuite une dérivée quelconque de l'une des fonctions inconnues, et désignant par $q$ un terme pris à volonté dans la suite $1, 2, \ldots, p$, nommons cote $q^{ième}$ de la dérivée en question l'entier obtenu en ajoutant à la cote $q^{ième}$ de la fonction inconnue les cotes homologues de toutes les variables de différentiation, distinctes ou non.
}
\end{quote}

\begin{quote}
\emph{
Denoting by
}

\medskip

\emph{
(1)\hspace{6cm} $x$, $y$, $\ldots$
}

\medskip

\emph{
the independent variables, and by 
}

\medskip

\emph{
(2)\hspace{6cm} $u$, $r$, $\ldots$
}

\medskip

\emph{
the unknown functions of any differential system, we make correspond to each of the quantities (1), (2)  $p$ integers, positive, zero or negative,  that we call, respectively, the first $\ll$\,cote\,$\gg$, second $\ll$\,cote\,$\gg$, $\ldots$, $p$-th $\ll$\,cote\,$\gg$ of this quantity. Then, considering any derivative of an unknown function, and denoting by $q$ a term taken freely from the sequence $1, 2, \ldots, p$, call the $q$-th $\ll$\,cote\,$\gg$ of the derivative in question the integer obtained by adding the homologous $\ll$\,cote\,$\gg$ of all of the variables, either distinct or not, of differentiation to the $q$-th $\ll$\,cote\,$\gg$ of the unknown function. 
}
\end{quote}

However, a complete algebraic formalization of this notion of cote wasn't obtained until 1929 by Janet in {\cite[\textsection 40]{Janet29}}, which we will recall in Section~\ref{SS:ParametricPrincipalJanet}. Moreover, he integrated the  notions of principal and parametric derivatives into a more general theory of orders on sets of derivatives, {\cite[Chapter II]{Janet29}}. The definitions for monomial orders given by Janet clarified the same notion previously introduced by Ch. Riquier in~\cite{Riquier93}. In particular, Janet made the notion of parametric and principal derivatives more explicit in order to distinguish the leading derivative in a polynomial PDE. In this way, he extended his algorithms on monomial PDE systems to the case of polynomial PDE systems. In particular, using these notions, he defined the property of completeness for a polynomial PDE system. Namely, a polynomial PDE system is complete if the associated set of monomials corresponding to leading derivatives of the system is complete. Moreover, he also extended the notion of complementary monomials to define the notion of \emph{initial conditions} for a polynomial PDE system as in the monomial case.

Finally, let us mention that Ch.~Riquier summarized known results on PDE systems in several books:~\cite{Riquier10} for PDE systems, \cite{Riquier28} for techniques of estimation.

\subsection{A finiteness result} 
In 1894, Arthur Tresse showed, as a preliminary result in the article \cite{Tresse94} on differential invariant theory, that PDE systems can be always reduced to systems of finitely many PDE. This is the first finiteness result relating to a module over a ring of differential operators. In particular, he showed in {\cite[Chap. I, Thm I]{Tresse94}} the following finiteness result:

\begin{quote}
\emph{
Un système d'équations aux dérivées partielles étant défini d'une manière quelconque, ce système est nécessairement limité, c'est-à-dire qu'il existe un ordre fini $s$, tel que, toutes les équations d'ordre supérieur à $s$ que comprend le système, se déduisent par de simples différentiations des équations d'ordre égal ou inférieur à $s$.
}
\end{quote}

\begin{quote}
\emph{
As a system of PDF might be defined arbitrarily, this system is necessarily limited, i.e., there exists a finite order, say $s$, such that all of the equations of order more than $s$ in the system can be deduced from simple differentiations of the equations whose order is less than or equal to $s$. 
}
\end{quote}

\subsection{Toward a more general normal form for PDE}
\label{SS:TowardMoreGeneralFormOfPDE}

Using the finiteness result of A. Tresse,   in \cite{Delassus1896} \'{E}. Delassus formalized and simplified,, Riquier's theory. In these works, one already finds an algorithmic approach to analysing ideals of the ring $\K[\frac{\partial}{\partial x_1}, \ldots, \frac{\partial}{\partial x_n}]$.
\'E. Delassus wrote {\cite[pp. 422-423]{Delassus1896}}:

\begin{quote}
\emph{
La solution du problème dépend de la recherche d'une forme canonique générale. M. Riquier, en faisant correspondre aux variables et aux inconnues des nombres entiers qu'il appelle cotes premières, cotes secondes, etc., est conduit à définir des systèmes orthonomes qu'il prend pour base de tous ses raisonnements. Il montre que tout système d'équations aux dérivées partielles peut se ramener à un système orthonome passif linéaire et du premier ordre. Dans de tels systèmes,
la formation par différentiation de toutes les équations, jusque à l'ordre infini, permet de séparer les dérivées des fonctions inconnues en deux classes, les unes étant principales et les autres paramétriques, et M. Riquier montre qu'en se donnant arbitrairement les valeurs initiales des dérivées paramétriques, on peut reconstruire les développements en séries des intégrales cherchées et que ces développements sont convergents.}

\emph{
Ces résultats sont établis en toute rigueur par M. Riquier, mais la démonstration, qu'il en donne, non seulement est très compliquée, mais est bien artificielle à cause de l'introduction de ces cotes qui interviennent d'une façon bien bizarre dans la question. Ceci justifierait déjà la publication de ce Travail où les résultats de M. Riquier sont retrouvés d'une façon beaucoup plus naturelle et plus simple en suivant une voie tout à fait différente; mais il y a plus, c'est que le Mémoire de M. Riquier n'a pas résolu la question aussi complètement qu'il est possible de le faire.
}
\end{quote}

\begin{quote}
\emph{The solution of the problem depends on how to find a general canonical form. By making correspond to variables and to unknown functions integers called first $\ll$\,cote\,$\gg$, second $\ll$\,cote\,$\gg$ etc., Mr. Riquier is led to define the 
$\ll$\,système orthonome\,$\gg$ which he takes as the base of his arguments.  
He shows that any system of partial differential equations can be reduced to a first-ordered passive linear $\ll$\,système orthonome\,$\gg$. 
In such systems, adding differentiations of all of the equations, up to infinite order,
allows one to separate the derivatives of the unknown functions into two classes, the one being principal and the other parametric. Mr. Riquier shows that given any initial values to the parametric derivatives, one can reconstruct the (formal) series expansion of the integrals we are looking for and that such series are convergent. 
}

\emph{These results are established completely rigorously by Mr. Riquier, but the proof is not only very complicated but is quite artificial owing to the introduction of these $\ll$\,cote\,$\gg$ which play a quite strange role in the question. 
This already may justify the publication of this work where the results of Mr. Riquier are recovered in a much more simple and natural way following a totally different path. But there is more, that is, the article of Mr. Riquier does not resolve the question as completely as  can be done. 
}
\end{quote}

\noindent Ch. Riquier answered to \'E. Delassus in {\cite[pp. 424]{Riquier97}}:

\begin{quote}
\emph{
Je m'étonne d'avoir été aussi peu compris. Que M. Delassus, retrouvant les résultats que j'ai le premier obtenus, estime y être arrivé par une voie plus simple, c'est une croyance que je m'explique chez lui, bien que je ne la partage pas, et que ses démonstrations me paraissent tout aussi compliquées que les miennes. Libre encore à M. Delassus de trouver « bizarre » l'attribution de cotes entières aux variables et aux inconnues, bien que cette idée ne me semble pas, à moi, plus singulière que celle de les ranger, comme il le fait, dans un ordre déterminé. Mais lorsqu'il soutient, et c'est là le point important de sa critique, que je n'ai pas résolu la question d'une manière complète, et qu'il est impossible, en suivant ma méthode, d'apercevoir << comment on pourrait grouper les coefficients arbitraires des développements des intégrales pour former des fonctions arbitraires, en
nombre fini, ayant avec ces dernières des relations simples >>, je ne puis, sans protester, laisser passer de semblables affirmations.
}
\end{quote}

\begin{quote}
\emph{
I am surprised I was so little understood. 
Mr. Delassus, recovering the results that I was the first to obtain, believes that
he arrived at the results in a much simpler way, this is what I believe that I understand of him, even if I don't really think so and his proofs seem quite as complicated as mine. M. Delassus is free to find strange, the attribution of total $\ll$\,cote\,$\gg$ to the variables and unknown functions, even though this idea seems to me no more singular than to order them, as he does, in a fixed order. But when he supports (his theory), and this is the important point in his criticism, that I did not solved the question in a complete manner, and that, by following my method, it is impossible to see ``how can one group all of the arbitrary coefficients of the series expansion of the integrals to form any function, in a finite number (of steps), having simple relations with these coefficients'', I could not ignore similar affirmations without protesting. }
\end{quote}

Apart from works of Ch. Riquier and \'E. Delassus, there had not been significant progress on the computation of normal forms for linear PDE systems. However, several monographs appeared on the topic and had a great influence on the community in the beginning of 20th century: Forsyth \cite{Forsyth90}, Weber~\cite{Weber00}, \'E. Goursat \cite{Goursat22}, Ch. Riquier \cite{Riquier10}. The research of new methods to compute normal forms of linear PDE systems was taken up by Janet in the period 1920-1930. 

\section{Algebraisation of monomial PDE systems}

The computational approach to reach normal forms for linear PDE systems in the work of 
Ch. Riquier and \'E. Delassus was not complete. The thesis of Janet provides a major contribution to the algebraisation of the problem considered by Ch. Riquier and \'E. Delassus by introducing an algorithmic method to compute normal forms of linear PDE systems. The procedure is based on a computation on a family of monomials associated to the PDE system. Finiteness properties on the set of monomials guarantee the termination of the procedure. In this section, we recall these constructions introduce by Janet. We recall also the results known by Janet on finiteness properties on set of monomials. In Section~\ref{SS:CompleteHigherOrderFiniteLineraPDESystems}, we will show how the results on monomials can be used to treat the general case of linear PDE systems.

\subsection{Monomial partial differential equation systems}
\label{SS:MonomialPartialDifferentialEquationsSystems}

In his thesis \cite{Janet20PhD}, Janet considered \emph{monomial PDE}, that is PDE of the following form
\begin{equation}
\label{Equation:PDEform}
\frac{\partial^{\alpha_1+\alpha_2+\ldots +\alpha_n}\varphi}{\partial x_1^{\alpha_1}\partial x_2^{\alpha_2}\ldots \partial x_n^{\alpha_n}}
\: = \:
f_{\alpha_1\alpha_2\ldots \alpha_n}(x_1,x_2,\ldots, x_n),
\end{equation}
where $\varphi$ is an unknown function and the $f_{\alpha_1\alpha_2\ldots \alpha_n}$ are analytic functions in several variables.
His objective was to compute an analytic function $\varphi$ which is a solution of the system.
He considered this problem, using an original algebraic approach, by seeing the differentiation operation as a multiplication operation on monomials. Tacitly, he used the ring isomorphism from the ring of polynomials in several variables with coefficients in an arbitrary field $\K$ to the ring of differential operators with constant coefficients. Note that, this isomorphism was established explicitly more than fifteen years later by W. Gröbner in {\cite[pp. 128]{Grobner37}} in a modern algebraic language as follows:

\begin{quote}
\emph{
Jedem Polynom $p(x)\subset \mathfrak{P}_n$ ordnen wir eineindeutig einen Differentialoperator $p\left(\frac{\partial}{\partial x}\right)$ zu, indem wir einfach die einzelnen Potenzprodukte $x_1^{i_1}\cdots x_n^{i_n}$ in $p(x)$ durch die Symbole 
$\frac{\partial^i}{\partial x_1^{i_1} \cdots \partial x_n^{i_n}}$, $(i=i_1+i_2+\cdots +i_n)$
ersetzen, was kurz durch $p(x) \leftrightarrow p\left(\frac{\partial}{\partial x}\right)$ angedeutet sei.
\\
Ist au$\beta$erdem auch $q(x)\leftrightarrow q\left(\frac{\partial}{\partial x}\right)$, so folgt leicht
\begin{align*}
p(x)+q(x) &\leftrightarrow p\left(\frac{\partial}{\partial x}\right) + q\left(\frac{\partial}{\partial x}\right)
\\
p(x)\cdot q(x) &\leftrightarrow p\left(\frac{\partial}{\partial x}\right) \cdot q\left(\frac{\partial}{\partial x}\right)
\end{align*}
Da bei dieser Zuordnung der Grundkörper $K$ elementweise festbleibt, unterscheiden sich die beiden Bereiche $\mathfrak{P}_n=K[x_1,\cdots, x_n]$ und 
$\mathfrak{D}_n=K\left[\frac{\partial}{\partial x_1},\cdots,\frac{\partial}{\partial x_n}\right]$ nur durch die verschiedene Bezeichnung ihrer transzendenten Elemente, sind also isomorph.
}
\end{quote}

\begin{quote}
\emph{
We assign each polynomial $p(x)\subset \mathfrak{P}_n$ to a differential operator $p\left(\frac{\partial}{\partial x}\right)$, by simply replacing a monomial $x_1^{i_1}\cdots x_n^{i_n}$ (appearing) in $p(x)$ with the symbol $\frac{\partial^i}{\partial x_1^{i_1} \cdots \partial x_n^{i_n}}$, $(i=i_1+i_2+\cdots +i_n)$, which is shortly expressed as 
 $p(x) \leftrightarrow p\left(\frac{\partial}{\partial x}\right)$.
If there is also $q(x)\leftrightarrow q\left(\frac{\partial}{\partial x}\right)$, it follows easily
that 
\begin{align*}
p(x)+q(x) &\leftrightarrow p\left(\frac{\partial}{\partial x}\right) + q\left(\frac{\partial}{\partial x}\right)
\\
p(x)\cdot q(x) &\leftrightarrow p\left(\frac{\partial}{\partial x}\right) \cdot q\left(\frac{\partial}{\partial x}\right)
\end{align*}
Since this assignment fixes the ground field $K$,  the two sets 
$\mathfrak{P}_n=K[x_1,\cdots, x_n]$ and
$\mathfrak{D}_n=K\left[\frac{\partial}{\partial x_1},\cdots,\frac{\partial}{\partial x_n}\right]$
differ only in the different names of their transcendental elements, they are isomorphic. 
}
\end{quote}

In this article, we will denote by 
\[
\Phi : \mathbb{K}[x_1,\ldots,x_n] \longrightarrow \mathbb{K}\left[\frac{\partial}{\partial x_1},\ldots,\frac{\partial}{\partial x_n}\right],
\]
the aforementioned ring isomorphism given explicitely by W. Gröbner from the ring of polynomials with $n$-variables to the ring of differential operators with constant coefficients.
Janet considered monomials in the variables $x_1,\ldots,x_n$ and implicitly used the isomorphism $\Phi$. 
In this way, he associated a monomial $x_1^{\alpha_1}x_2^{\alpha_2}\ldots x_n^{\alpha_n}$ to the differential operator  
\[
\frac{\partial^{\alpha_1+\alpha_2+\ldots +\alpha_n}\quad}{\partial x_1^{\alpha_1}\partial x_2^{\alpha_2}\ldots \partial x_n^{\alpha_n}}.
\]
In his thesis~{\cite[Chapitre I]{Janet20PhD}}, Janet considered \emph{monomial PDE systems}, that is those whose equations are of the form~(\ref{Equation:PDEform}), and which have finitely many such equations.
Such a system can be written as the following family:
\begin{equation}
\label{Equation:MonomialPDESystem}
(\Sigma)
\qquad
\frac{\partial^{\alpha_1+\alpha_2+\ldots +\alpha_n}\;\varphi}{\partial x_1^{\alpha_1}\partial x_2^{\alpha_2}\ldots \partial x_n^{\alpha_n}}
\: = \:
f_{\alpha_1,\ldots,\alpha_n}(x_1,x_2,\ldots, x_n), \qquad (\alpha_1,\ldots,\alpha_n) \in I, 
\end{equation}
where $\varphi$ is an unknown function and the $f_{\alpha_1,\ldots,\alpha_n}$ are analytic functions in several variables, and indexed by a finite subset $I$ of $\Nb^n$. 

\subsection{Finiteness properties on monomials}
\label{SS:FinitenessPropertiesMonomials}

Using the ring isomorphism $\Phi$ defined above, Janet associated a PDE system $(\Sigma)$ of the form (\ref{Equation:MonomialPDESystem}) to the set $\lm(\Sigma)$ of monomials defined as follows
\[
\lm(\Sigma) \: = \: \{\;x_1^{\alpha_1}\ldots x_n^{\alpha_n}\;\; |\;\;(\alpha_1,\ldots,\alpha_n)\in I\}.
\]
In his hypotheses, Janet excluded the case in which the system has an infinite number of equations. Indeed, there are finiteness results that he stated as the \emph{Théorème général sur certaines suites de monomes}, {\cite[\textsection 1]{Janet20}}:
\begin{quote}
\emph{Une suite des monomes $M_1, M_2, \ldots$ telle que chacun d'entre eux n'est multiple d'aucun des précédents ne comprend qu'un nombre fini de monomes.
}
\end{quote}

\begin{quote}
\emph{A sequence of monomials $M_1$, $M_2$, $\ldots$ such that each monomial is not a multiple of any preceding one contains only a finite number of monomials.}
\end{quote}

He proved this theorem by induction on the number of variables constituting the monomials.
Janet considered these finiteness properties with the objective of giving an inductive form to his constructions.
Note that the finiteness result on PDE systems was already published in 1894 by Tresse in \cite{Tresse94}, and used by \'E. Delassus as exposed in Section~\ref{SS:TowardMoreGeneralFormOfPDE}.
However, the finiteness assumption in this context was formulated algebraically for the first time by Janet. This result had already been known by Leonard Eugene Dickson  in {\cite[Lemma A]{Dickson13}}. 

\begin{quote}
\emph{
Lemma A. Any set $S$ of functions of the type
\[
F \: = \: x_1^{e_1}x_2^{e_2}\ldots x_n^{e_n},
\qquad
\text{\emph{(e's integers $\geqq 0$)}}
\qquad\qquad (1)
\]
contains a finite number of functions $F_1,\ldots,F_k$ such that each function $F$ of the set $S$ can be expressed as a product $F_if$, where $f$ is of the form (1), but is not necessarily in the set $S$.
}
\end{quote}
This result was published in 1913 in an article on number theory in the American journal \emph{American Journal of Mathematics}, but due to the First World War, it would take a long time before these works were accessible to French mathematical community.

The results presented by Janet in his thesis follows those of Ch. Riquier, with an original algebraic formulation. The new algebraic approach to this well-studied problem in PDE systems proposed by Janet was made possible by the influence of the German mathematical school on the academic development of Janet. In the introduction of his thesis, {\cite[Introduction \textsection 2]{Janet20PhD}}, he presented his contribution as follows:
\begin{quote}
\emph{
Le présent travail a pour objet essentiel l'exposition simple des résultats de M. Riquier. Cette exposition nous conduira naturellement à certains résultats de nature algébrique qui complètent la théorie des formes donnée par M. Hilbert.
}
\end{quote}

\begin{quote}
\emph{The main purpose of this work is a simple presentation of the results of Mr. Riquier. This exposition leads us naturally to certain results of algebraic nature which complete the theory of polynomials given by Mr. Hilbert.
}
\end{quote}

Here, Janet mentions the finiteness result of D. Hilbert on what we today call the Noetherian character of the polynomial ring over a Noetherian ring, now called \emph{Hilbert's basis theorem}, and published in \cite{Hilbert1890}.

\subsection{On algebraic finiteness properties}
\label{SS:OnAlgebraicFinitenessProperties}

The constructions of Janet are based on some remarkable properties on monomial ideals that he developed in his thesis, \cite{Janet20PhD}, and published in \cite{Janet20} and \cite{Janet20a}.
In particular, as explained above, he gave another formulation of Dickson's Lemma on the finiteness of generating sets of monomial ideals.
This finiteness property is essential for Noetherian properties on the set of monomials. Note that Janet was not familiar with the axiomatisation of the algebraic structure of ideal and the property of Noetherianity introduced by E. Noether at the same time in \cite{Noether1921} and \cite{Noether23}.

The finiteness property \emph{Théorème général sur certaines suites de monomes} recalled above, was formulated by Janet by introducing the property, of a family of monomials $\Ur$, of being \emph{multiplicatively stable}, which means that $\Ur$ is closed under multiplication by monomials in $\Mr(x_1,\ldots,x_n)$.
By this finiteness property,  if $\Ur$ is a multiplicatively stable, then it contains only finitely many elements which are not multiples of any other elements in $\Ur$. Hence, there exists a finite subset $\Ur_f$ of $\Ur$ such that for any $u$ in $\Ur$, there exists $u_f$ in $\Ur_f$ such that $u_f$ divides $u$.
From the finiteness property, Janet deduced the \emph{ascending chain condition} on multiplicatively stable monomial sets that he formulated as follows. Any ascending sequence of multiplicatively stable subsets of $\Mr(x_1,\ldots,x_n)$
\[
\Ur_1 \subset \Ur_2 \subset \; \ldots \; \subset \Ur_k \subset \ldots 
\]
is finite. This corresponds to the \emph{Noetherian property} introduced by E. Noether in {\cite[\textsection 1]{Noether1921}} in the following terms
\begin{quote}
\emph{
Satz I (Satz von der endlichen Kette): Ist $\mathfrak{M}$, $\mathfrak{M}_1$, $\mathfrak{M}_2$, \ldots, $\mathfrak{M}_\nu$, \ldots ein abzählbar unendliches System von Idealen in $\Sigma$, von denen jedes durch das folgende teilbar ist, so sind von einem endlichen Index $n$ an alle Ideale identisch, $\mathfrak{M}_n=\mathfrak{M}_{n+1}=\ldots$ M. a.W.: Bildet  $\mathfrak{M}$, $\mathfrak{M}_1$, $\mathfrak{M}_2$, \ldots, $\mathfrak{M}_\nu$, \ldots eine einfach geordnete Kette von Idealen derart, da$\beta$ jedes Ideal ein echter Teiler des unmittelbar vorangehenden ist, so bricht die Kette im EndIichen ab.
}
\end{quote}

\begin{quote}
\emph{
Theorem I (theorem of finite chain): Let $\mathfrak{M}$, $\mathfrak{M}_1$, $\mathfrak{M}_2$, \ldots, $\mathfrak{M}_\nu$, \ldots be a system of countably infinite ideals in $\Sigma$, each of which is divisible by the next ideal. Then, there exists a finite index $n$ from which all of the ideals are identical, $\mathfrak{M}_n=\mathfrak{M}_{n+1}=\ldots$ In other words, let us form $\mathfrak{M}$, $\mathfrak{M}_1$, $\mathfrak{M}_2$, \ldots, $\mathfrak{M}_\nu$, \ldots a simply ordered chain of ideals as above so that each ideal is divisible by the next ideal, then the chain stops after a finite number of steps. 
}
\end{quote}

\subsection{On the notion of module}
\label{SS:OnNotionOfModule}
Throughout his work on the analysis of PDE and until his monograph~\cite{Janet29} appeared in 1929, Janet developed computational methods to deal with monomials and polynomials over a field. Nowadays, these methods are known and developed in the language of ideals. The use of a formal definition of the notion of ideal appeared progressively in the series of Janet's works on formal analysis of linear PDE systems, \cite{Janet20PhD, Janet21a, Janet21b}.
Note that, at this time, Janet knew only the structure of ideal of the ring of integers of number field.
The first formulation of the structure of ideal appeared in the series of articles by Richard Dedekind, \cite{Dedekind77}, see also {\cite[\textsection 177]{LejeuneDirichlet94}}. 
Hilbert investigated in a systematic way the notion of ideal of a ring of commutative polynomials of several variables in a seminal paper \cite{Hilbert1890} under the terminology of \emph{algebraic forms}. In particular, he proved such results as the ring of polynomials over a field is Noetherian, now called Hilbert's basis theorem. 
Notice that N. M. Gunther dealt with such a structure in \cite{Gunther13b}.
The modern algebraic formulation of the notion of ideal over a general commutative ring was only introduced in 1921 by E. Noether in \cite{Noether1921}. 

In the case of monomial PDE systems, Janet explained his constructions without using the structure of monomial ideal in the sense of an ideal generated by monomials. Instead, his results are formulated using the notion of multiplicative cone.
In his thesis, {\cite[Chapter I, \textsection 3]{Janet20}}, Janet defined the notion of \emph{module de monomes} (module of monomials) by specifying its finiteness properties.

\begin{quote}
\emph{
Nous dirons qu'un système de monomes forme un module si tout multiple d'un de ces monomes appartient au système. Un module est toujours constitué par les multiples d'un nombre fini de monomes. Nous dirons quelquefois que ces monomes forment une base pour le module.
}
\end{quote}

\begin{quote}
\emph{We say that a system of monomials constitutes a $\ll$ module $\gg$ if any multiple of one of these monomials belongs to the system. A module always consists of the multiples of a finite number of monomials. We sometimes say that these monomials form a base of the module.
}
\end{quote}
In this note, \emph{module de monomes} will be called multiplicative cone, and this notion will be presented in the next section. 

In an article published in 1924, \cite{Janet24}, Janet used the notion of \emph{algebraic form}, introduced by Hilbert, in his study of linear polynomial PDE systems, that is a PDE system where each equation is defined by a polynomial in partial differential operators.
In this polynomial situation, he used the structure of polynomial ideal as D. Hilbert did.
Indeed, following the approach developed by D. Hilbert in {\cite[IV. \emph{Die charakteristische Function eines Moduls}]{Hilbert1890}}, Janet recalled in {\cite[Chapter III, \textsection 52]{Janet29}} the definition of \emph{polynome caractéristique ou la postulation} of the module of forms of a polynomial PDE system. In modern language, this polynomial corresponds to the coefficients of the Hilbert series of the ideal generated by a polynomial PDE system.
He used such Hilbert series to define the property of involutivity on polynomial PDE systems in Chapter III of his monograph.
In addition to his work on the solvability of linear PDE systems, in a series of publications \cite{Janet13a,Janet22a,Janet24}, Janet studied the notion of character and involutivity of linear PDE systems. 
We do not develop the results obtained by Janet in this direction.

\section{Janet's completion procedure}
\label{Section:JanetCompletion}

We present the main algorithmic ingredient in the construction of Janet, namely the completion procedure of a set of monomials with respect the notion of multiplicative variable. The completeness property is formulated using the notion of multiplicative cone, and thus can be characterized using the notion of involutive division. In this section, we recall these constructions of Janet on a set of monomials, which were mainly introduced in the memoir of his thesis.

\subsection{Multiplicative cone of a set of monomials}
\label{SS:MultiplicativeCone}
For a finite set $\Ur$ of monomials in variables $x_1,\ldots, x_n$, Janet gave an inductive construction of the multiplicative cone $\cone(\Ur)$ generated by~$\Ur$, that is the set of monomials $u$ such that there exists $u'$ in $\Ur$ that divides $u$.
With the objective of introducing the involutive cone of a set of monomials as a refinement of the multiplicative cone, Janet gave an inductive construction of $\cone(\Ur)$ as follows. First, he defined, for every $0\leq \alpha_n \leq \deg_n(\Ur)$,
\[
[\alpha_n] \: = \: \{u \in \Ur\;|\; \deg_{n}(u) = \alpha_n\,\},
\]
in such a way, that the family $([0],\ldots ,[\deg_n(\Ur)])$ forms a partition of $\Ur$.
By setting, for every $0\leq \alpha_n \leq \deg_n(\Ur)$,
\[
\overline{[\alpha_n]} \: = \:\{u \in \Mr(x_1,\ldots,x_{n-1}) \;|\; ux_n^{\alpha_n} \in \Ur\,\},
\]
he defined for every $0\leq i \leq \deg_n(\Ur)$
\[
\Ur_i' \: = \: \underset{0\leq \alpha_n \leq i}{\bigcup} \{u \in \Mr(x_1,\ldots,x_{n-1}) \;|\; \text{there exists $u'\in \overline{[\alpha_n]}$ such that $u'|ux_n^{\alpha_n}$}\,\}.
\]
By denoting 
\[
\Ur_k \: = \:
\begin{cases} 
\{\,ux_n^k\;|\; u\in \Ur'_k\,\} & \text{if $k< \deg_n(\Ur)$,} \\ 
\{\,ux_n^k\;|\; u\in \Ur'_{\deg_n(\Ur)}\,\} & \mbox{if $k\geq \deg_n(\Ur)$.}
\end{cases} 
\]
he constructed the multiplicative cone $\cone(\Ur)$ as the set $\underset{k\geq 0}{\bigcup} \Ur_k$. 

\subsection{The notion of multiplicative variable}
\label{SS:NotionMultiplicativeVariable}
In 1920, Janet introduced the notion of \emph{multiplicative variable}, see~{\cite[\textsection 7]{Janet20}} and~{\cite[\textsection 1]{Janet20a}}. In \cite{Janet20a}, he wrote

\begin{quote}
\emph{
Soit un système formé d'un nombre fini de \emph{monomes} $(M)$ à $n$ variables $x_1,x_2,\ldots,x_n$ ; $x_i$ sera dite \emph{multiplicatrice} pour $\overline{M}=x_n^{\alpha_n}x_{n-1}^{\alpha_{n-1}}\ldots x_1^{\alpha_1}$ dans le système $(M)$ si parmi les $(M)$ où $x_n,x_{n-1},\ldots,x_{i+1}$ ont les exposants $\alpha_n,\alpha_{n-1},\ldots,\alpha_{i+1}$, il n'y en a pas où $x_i$ ait un exposant supérieur à $\alpha_i$; on dira qu'un monome provient de $\overline{M}$ s'il est le produit de $\overline{M}$ par un monome ne contenant que des variables multiplicatrices de $\overline{M}$. 
}
\end{quote}

\begin{quote}
\emph{
Let $(M)$ be a system made of a finite number of \emph{monomials} of $n$ variables
$x_1, x_2, \ldots, x_n$ ; $x_i$ will be \emph{multiplicative} for  $\overline{M}=x_n^{\alpha_n}x_{n-1}^{\alpha_{n-1}}\ldots x_1^{\alpha_1}$ in the system $(M)$ if
among $(M)$ where the exponents of $x_n,x_{n-1},\ldots,x_{i+1}$ are $\alpha_n,\alpha_{n-1},\ldots,\alpha_{i+1}$ there is no monomial where $x_i$ has an exponent  greater than $\alpha_i$; we say that a monomial comes from $\overline{M}$
if it is the product of $\overline{M}$ with a monomial which contains only multiplicative variables of $\overline{M}$.
}
\end{quote}

This definition can be expanded as follows. 
Given a finite set $\Ur$ of monomials in the variables $x_1,\ldots,x_n$, we define,
for all $1\leq i \leq n$, the following subset of $\Ur$:
\[
[\alpha_i,\ldots,\alpha_n] \: = \: 
\{u\in \Ur\;|\;\deg_j(u)=\alpha_j\;\;\text{for all}\;\;i\leq j \leq n\}.
\]
That is $[\alpha_i,\ldots,\alpha_n]$ contains monomials of $\Ur$ of the form $vx_i^{\alpha_i}\ldots x_n^{\alpha_n}$, with $v$ in $\Mr(x_1,\ldots,x_{i-1})$.
The sets $[\alpha_i,\ldots,\alpha_n]$, for $\alpha_i,\ldots,\alpha_n$ in $\Nb$, form a partition of $\Ur$. Moreover, for all $1\leq i \leq n-1$, we have $[\alpha_i,\alpha_{i+1},\ldots,\alpha_n] \subseteq [\alpha_{i+1},\ldots,\alpha_n]$ and the sets $[\alpha_i,\ldots,\alpha_n]$, where $\alpha_i\in \Nb$, form a partition of~$[\alpha_{i+1},\ldots,\alpha_n]$.

The variable $x_n$ is said to be \emph{multiplicative} for a monomial $u$ in $\Ur$, if $\deg_n(u) = \deg_n(\Ur)$.
For $i\leq n-1$, the variable $x_i$ is said to be \emph{multiplicative} for $u$ if 
\[
u\in [\alpha_{i+1},\ldots,\alpha_n]
\qquad\text{and}\qquad
\deg_i(u)=
\deg_i([\alpha_{i+1},\ldots,\alpha_n]).
\]
We will denote by $\mult_{\div{J}}^\Ur(u)$ the set of multiplicative variables of $u$ with respect to the set~$\Ur$. 
The set of \emph{non-multiplicative variables of $u$ with respect to the set~$\Ur$}, denoted by $\nonmult_{\div{J}}^{\Ur}(u)$, is defined as the complementary set of $\mult_{\div{J}}^\Ur(u)$ in the set $\{x_1,\ldots,x_n\}$.

The notion of multiplicative variable is local in the sense that it is defined with respect to a subset~$\Ur$ of the set of all monomials. A monomial $u$ in $\Ur$ is said to be a \emph{Janet divisor} of a monomial $w$ with respect to~$\Ur$, if $w=uv$ and all variables occurring in $v$ are multiplicative with respect to $\Ur$. In this way, we distinguish the set $\cone_{\div{J}}(\Ur)$ of monomials having a Janet divisor in~$\Ur$, called \emph{$\div{J}$-multiplicative} or \emph{involutive cone} of $\Ur$, from the set $\cone(\Ur)$ of multiple of monomials in $\Ur$ for the classical division.
Explicitly, the involutive cone is defined by 
\[
\cone_{\div{J}}(\Ur) \: = \: 
\underset{u \in \Ur}{\bigcup} \; \{\, uv \; | \; v\in \Mr(\mult_{\div{J}}^\Ur(u))\,\}.
\]

\subsection{Completeness of a set of monomials}
\label{SS:CompletenessSetMonomials}
Janet introduced, in~{\cite[\textsection 1]{Janet20a}}, the notion of \emph{completeness} of a set of monomials:
\begin{quote}
\emph{
Un monome ne peut provenir de deux monomes $(M)$ différents. 
Pour que tout multiple d'un monome du système provienne d'un de ces monomes, il faut et il suffit qu'il en soit ainsi de tous les produits obtenus en multipliant un $(M)$ par une de ses variables non-multiplicatrices.
Lorsque cette condition sera réalisée, le système $(M)$ sera dit complet.
}
\end{quote}

\begin{quote}
\emph{
A monomial cannot come  from two different monomials in 
$(M)$. Any multiple of a monomial in the system is deduced from one of these monomials if and only if any product of monomials obtained by multiplying a monomial in $(M)$ with one of its non-multiplicative variables is deduced from a monomial in $(M)$. When this condition is realized, the system $(M)$ is said to be complete.
}
\end{quote}

In this formulation, the meaning of \emph{provenir} (come from) can be explained as follows. A monomial $v$ comes from a monomial $u$ if $v$ can be decomposed into a product $v=uw$, where all the variables in $w$ belong to $\mult_{\div{J}}^\Ur(u)$.
In the above formulation of completeness, the notion of involutive cone of a set of monomials $\Ur$ appears implicitly. Janet division being a refinement of the classical division, the set $\cone_{\div{J}}(\Ur)$ is a subset of $\cone(\Ur)$. Janet called a set of monomials $\Ur$ \emph{complete} precisely when this inclusion is an equality, namely when the involutive cone is equal to the set of all products $uv$ of monomials such that $u$ is in $\Ur$ and $v$ is an arbitrary monomial.  
He thus obtained a characterization of completeness of a finite set of monomials. He proved, cf. {\cite[pp. 20]{Janet29}}, that a finite set $\Ur$ of monomials is complete if and only if, for any $u$ in $\Ur$ and any non-multiplicative variable $x$ of $u$ with respect to $\Ur$, $ux$ is in $\cone_{\div{J}}(\Ur)$.

Using this characterization, Janet deduced in {\cite[pp. 21]{Janet29}} a \emph{completion procedure} for any finite set $\Ur$ of monomials in $\Mr(x_1,\ldots,x_n)$, whose principle consists in adding monomials $ux$, for all $u$ in $\Ur$ and $x\in{\nonmult_{\div{J}}^{\Ur}(u)}$, such that $ux$ is not in $\cone_{\div{J}}(\Ur)$ and iterating this process until the set contains no such~$ux$ with this property.

With this constructive approach, he proved, cf. {\cite[pp. 21]{Janet29}}, that for any finite set $\Ur$ of monomials there exists a finite complete set $J(\Ur)$ that contains $\Ur$ and such $\cone(\Ur) = \cone(J(\Ur))$.
Note that Janet does not give a proof of the termination of the completion procedure. 

In order to illustrate this construction, let us recall an example from {\cite[pp. 28]{Janet29}}.
Consider $\Ur=\{\,x_3x_2^2,x_3^3x_1^2\,\}$. The following table gives the multiplicative variables for the monomials of $\Ur$:
\begin{center}
\begin{tabular}{c|ccc}
$x_3^3x_1^2$ & $x_3$ & $x_2$ & $x_1$\\
$x_3x_2^2$ &  & $x_2$ & $x_1$\\
\end{tabular}
\end{center}
The set $\Ur$ can then be completed as follows. The monomial $(x_3x_2^2)x_3$ is not in $\cone_\div{J}(\Ur)$; we set $\widetilde{\Ur}\leftarrow \Ur\cup\{x_3^2x_2^2\}$ and we compute the multiplicative variables with respect to $\widetilde{\Ur}$:
\begin{center}
\begin{tabular}{c|ccc}
$x_3^3x_1^2$ & $x_3$ & $x_2$ & $x_1$\\
$x_3^2x_2^2$ &  & $x_2$ & $x_1$\\
$x_3x_2^2$ &  & $x_2$ & $x_1$\\
\end{tabular}
\end{center}
The monomial $(x_3x_2^2)x_3$ is in $\cone_\div{J}(\widetilde{\Ur})$, but $(x_3^2x_2^2)x_3$ is not in $\cone_\div{J}(\widetilde{\Ur})$;  we set 
$\widetilde{\Ur} \leftarrow \widetilde{\Ur} \cup \{x_3^3x_2^2\}$. The multiplicative variables of this new set of monomials are
\begin{center}
\begin{tabular}{c|ccc}
$x_3^3x_2^2$ & $x_3$ & $x_2$ & $x_1$\\
$x_3^3x_1^2$ & $x_3$ &  & $x_1$\\
$x_3^2x_2^2$ &  & $x_2$ & $x_1$\\
$x_3x_2^2$ &  & $x_2$ & $x_1$\\
\end{tabular}
\end{center}
The monomial $(x_3^3x_1^2)x_2$ is not in $\cone_\div{J}(\widetilde{\Ur})$, the other products are in $\cone_\div{J}(\widetilde{\Ur})$, and we prove that the system
\[
\widetilde{\Ur} \: = \:
\{\,x_3^3x_1^2,x_3x_2^2,x_3^2x_2^2,x_3^3x_2^2,x_3^3x_2x_1^2\,\}
\]
is complete, so $J(\Ur) = \widetilde{\Ur}$.

\section{Initial value problem}

Given an ideal generated by a set of monomials, Janet distinguished the family of monomials contained in the ideal and those contained in the complement of the ideal. The notion of multiplicative and non-multiplicative variables is used to stratify these two families of monomials. This leads to a refinement of the classical division on monomials. These constructions are based on the notion of complementary monomial defined as follows.

\subsection{Complementary monomials}
\label{SS:ComplementaryMonomials}
The notion of \emph{complementary monomial} appear for the first time in~{\cite[\textsection 1]{Janet20a}}. He wrote
\begin{quote}
\emph{
[...] étant donné un système quelconque de monomes (M), on est en possession d'un procédé régulier pour répartir respectivement : $I^o$ tous les monomes multiples d'un M au moins ; $2^o$ tous les autres monomes, en un nombre fini d'ensembles sans éléments communs, les monomes d'un ensemble se déduisant d'un monome déterminé en le multipliant par tous les monomes ne contenant que certaines variables déterminées.
}
\end{quote}

\begin{quote}
\emph{
[...] given any system of monomials (M), one has a regular procedure to divide respectively to: $I^o$ any multiple of at least one monomial in $M$ ; 
$2^o$ all of the other monomials, to a finite number of sets without common elements, 
the monomials of a set can be obtained from a given monomial by multiplying all of the monomials containing only specific variables.
}
\end{quote}

This notion was made explicit in {\cite[\textsection 2]{Janet21a}}.
The set of complementary monomials of a set of monomials~$\Ur$ is the set of monomials denoted by $\comp{\Ur}$ defined by the following disjoint union
\begin{equation}
\label{Equation:ComplementaryMonomials}
\comp{\Ur} \: = \: \bigcup_{1\leq i \leq n} \; \compp{i}{\Ur},
\end{equation}
where 
\[
\compp{n}{\Ur} \: = \: 
\{ x_n^\beta \; | \; 0 \leq \beta \leq \deg_n(\Ur)\;\text{and}\; 
[\beta]= \varnothing \},
\]
and for every $1\leq i < n$,
\begin{align*}
\compp{i}{\Ur} \: = \:
\big\{\,
x_i^\beta x_{i+1}^{\alpha_{i+1}}\ldots x_n^{\alpha_n} \; \big| \;
[\alpha_{i+1},\ldots,\alpha_n] \neq \varnothing, 
\; \hspace{6cm}
\\
0 \leq \beta < \deg_i([\alpha_{i+1},\ldots,\alpha_n]),
\;
[\beta, \alpha_{i+1}, \ldots , {\alpha_n}] = \varnothing
\,\big\}.
\end{align*}

For any monomial $u$ in $\comp{\Ur}$, we define the set $\cmult_{\div{J}}^{\comp{\Ur}}(u)$ of \emph{multiplicative variables for $u$ with respect to complementary monomials} in $\comp{\Ur}$ as follows. If the monomial $u$ is in $\compp{n}{\Ur}$, we set
\[
\cmult_{\div{J}}^{\compp{n}{\Ur}}(u) \: = \: \{x_1,\ldots,x_{n-1}\}.
\]
For $1\leq i \leq n-1$, for any monomial $u$ in $\compp{i}{\Ur}$, there exists $\alpha_{i+1},\ldots,\alpha_n$ such that $u\in [\alpha_{i+1},\ldots,\alpha_n]$. Then 
\[
\cmult_{\div{J}}^{\compp{i}{\Ur}}(u) \: = \: \{x_1,\ldots,x_{i-1}\} \cup 
\mult_{\div{J}}^\Ur([\alpha_{i+1},\ldots,\alpha_n]).
\] 
Finally, for $u$ in $\comp{\Ur}$, there exists an unique $1\leq i_u \leq n$ such that $u\in \compp{i_u}{\Ur}$. Then we set
\[
\cmult_{\div{J}}^{\comp{\Ur}}(u) \: = \: \cmult_{\div{J}}^{\compp{i_u}{\Ur}}(u).
\]
We define the involutive cone of the complementary family of a family $\Ur$ of monomials as follows
\[
{\cone}^{\complement}_{\div{J}}(\Ur) \: = \: 
\underset{u \in\comp{\Ur}}{\bigcup} \; \{\, uv \; | \; v\in \Mr(\cmult_{\div{J}}^{\comp{\Ur}}(u))\,\}.
\]
Janet proved, cf. {\cite[pp. 18]{Janet29}}, that for any finite set $\Ur$ of monomials of $\Mr(x_1,\ldots,x_n)$, we have the following partition 
\begin{equation}
\label{E:MacaulayMonomial}
\Mr(x_1,\ldots,x_n) \: = \: \cone(\Ur) \amalg {\cone}^{\complement}_{\div{J}}(\Ur).
\end{equation}
An other form of this equality in the case of polynomial ideals was proved by Francis Sowerby Macaulay in \cite{Macaulay27}.

\subsection{The space of initial conditions}
During the 1920s, Janet's works are  mainly concerned with the analysis of \emph{Cauchy's problems}. That is, the problem of proving the existence and the uniqueness of solutions for PDE systems under given initial conditions. In \cite{Janet21a, Janet21b} he considered the complete integrability problem of monomial PDE systems. In particular, in {\cite[pp. 244]{Janet21a}} he formulated the problem as follows:

\begin{quote}
\emph{
Proposons-nous de déterminer une fonction $u$ telle que celles de ses dérivées qui sont caractérisées par les monômes $(M)$ d'un système complet donné soient des fonctions données des $n$ variables indépendantes $x_1,x_2,\ldots x_n$. Nous apercevons immédiatement certaines conditions de possibilité du problème : à chacune des identités $\overline{M}.x_i= M.x_1^{\alpha_1}x_2^{\alpha_2}\ldots x_n^{\alpha_n}$ que mentionne la définition précédente correspond une relation entre les fonctions auxquelles on cherche à égaler les dérivées correspondant aux $(M)$:
\[
\frac{\partial\overline{f}}{\partial x_i}
\: = \:
\frac{\partial^{\alpha_1+\alpha_2+\ldots +\alpha_n}\;f}{\partial x_1^{\alpha_1}\partial x_2^{\alpha_2}\ldots \partial x_n^{\alpha_n}}
\qquad\qquad
\text{[conditions (I)]}
\]
(si du moins on suppose la continuité des dérivées de $u$ que fait intervenir l'égalité précédente).
}
\end{quote}

\begin{quote}
\emph{Let us propose to determine a function $u$ such that those of  its derivatives that are characterized by the monomials $(M)$ of a given complete system shall be the given functions of $n$ independent variables $x_1,x_2,\ldots x_n$. We see immediately certain conditions of possibilities of the problem: to each identity $\overline{M}.x_i= M.x_1^{\alpha_1}x_2^{\alpha_2}\ldots x_n^{\alpha_n}$ that mentions the precedent definition corresponds a relation between the functions that we are searching for 
to make the equality of the corresponding derivatives to $(M)$:
\[
\frac{\partial\overline{f}}{\partial x_i}
\: = \:
\frac{\partial^{\alpha_1+\alpha_2+\ldots +\alpha_n}\;f}{\partial x_1^{\alpha_1}\partial x_2^{\alpha_2}\ldots \partial x_n^{\alpha_n}}
\qquad\qquad
\text{[conditions (I)]}
\]
(if at least we suppose the continuity of the derivatives of $u$ that appear in above equality).
}
\end{quote}

In \cite{Janet21b}, Janet considered monomial PDE systems of the form (\ref{Equation:MonomialPDESystem}), which he supposed to be finite using the arguments presented in Section~\ref{SS:FinitenessPropertiesMonomials}.
In Section \ref{SS:MonomialPartialDifferentialEquationsSystems}, we recalled the way in which Janet associated to each monomial $x^\alpha$ in variables $x_1,\ldots,x_n$ a differential operator $D^\alpha$ via the isomorphism~$\Phi$. 
In this way, to a monomial PDE system $(\Sigma)$ on variables $x_1,\ldots,x_n$ he associated a finite set $\lm(\Sigma)$ of monomials. Using the completion procedure recalled in Section~\ref{SS:CompletenessSetMonomials}, he showed that any such set $\lm(\Sigma)$ of monomials can be completed into a finite complete set $J(\lm(\Sigma))$ having the same multiplicative cone as $\lm(\Sigma)$.

Suppose that the set of monomials $\lm(\Sigma)$ is finite and complete. We have 
\[
\cone(\lm(\Sigma)) \: = \: \cone_{\div{J}}(\lm(\Sigma))
\]
Thus, for any monomial~$u$ of $\lm(\Sigma)$ and non-multiplicative variable $x_i$ in $\nonmult_\div{J}^{\lm(\Sigma)}(u)$, there exists a decomposition
\[
ux_i = vw,
\]
where $v$ is in $\lm(\Sigma)$ and $w$ belongs to $\Mr(\mult_{\div{J}}^{\lm(\Sigma)}(v))$.
For any such decomposition, it corresponds to a \emph{compatibility condition} of the monomial PDE system $(\Sigma)$, that is,  for $u=x^{\alpha}$, $v=x^{\beta}$ and $w=x^{\gamma}$
with $\alpha, \beta$ and $\gamma$ in $\Nb^n$, 
\begin{equation}
\label{E:ConditionI}
\frac{\partial f_{\alpha}}{\partial x_i} 
\: = \:
D^\gamma f_{\beta}.
\end{equation}
This condition corresponds to the \emph{conditions (I)} above mentioned by Janet.
Let us denote by $(C_\Sigma)$ the set of all such compatibility conditions.
Janet showed that with the completeness hypothesis, this set of compatibility conditions is sufficient for the monomial PDE system $(\Sigma)$ to be integrable. 

Let us consider the set $\comp{\lm(\Sigma)}$ of complementary monomials of the finite complete set~$\lm(\Sigma)$. Suppose that the monomial PDE system $(\Sigma)$ satisfies all the compatibility conditions in $(C_\Sigma)$. Under this hypothesis, Janet associated to each monomial $v=x_1^{\beta_1}\ldots x_n^{\beta_n}$ of $\comp{\lm(\Sigma)}$ an analytic function 
\[
\varphi_{\beta_1,\ldots,\beta_n}(x_{i_1},\ldots,x_{i_{k_v}}),
\]
where $\{x_{i_1},\ldots,x_{i_{k_v}}\}=\cmult_{\div{J}}^{\comp{\lm(\Sigma)}}(v)$. 
As a consequence of the decomposition~(\ref{E:MacaulayMonomial}), the set of such analytic functions provides a compatible initial condition. In {\cite[\textsection 7]{Janet21a}},  he obtained the following solvability result:

\begin{quote}
\emph{
Supposons que ces conditions (I) soient réalisées. Si le problème posé a une solution, cette solution vérifie bien évidemment, en particulier, les équations obtenues en annulant dans chacune
des équations proposées les variables non multiplicatrices du
premier membre. Réciproquement, considérons une solution des équations ainsi obtenues, je dis qu'elle est solution des équations proposées.
}
\end{quote}

\begin{quote}
\emph{Suppose that these conditions (I) are realized. 
If the given problem has a solution, this solution verifies evidently, in particular, the equations obtained by eliminating, in each of the given equations, the non-multiplicative variables in the left-hand side. Conversely, consider a solution of thus obtained equations, I say that it is a solution of the given equations.  
}
\end{quote}

Using the notations above on complementary monomials, this result can be formulated as follows.

\begin{theorem}
\label{Theorem:BoundaryConditions}
Let $(\Sigma)$ be a finite monomial PDE system such that $\lm(\Sigma)$ is complete. If~$(\Sigma)$ satisfies the compatibility conditions $(C_\Sigma)$, then it always admits a unique solution with initial conditions given for any $v=x_1^{\beta_1}\ldots x_n^{\beta_n}$ in $\comp{\lm(\Sigma)}$ by 
\[
\left. 
\frac{\partial^{\beta_1+\beta_2+\ldots +\beta_n}\;\varphi}{\partial x_1^{\beta_1}\partial x_2^{\beta_2}\ldots \partial x_n^{\beta_n}}
\right|_{x_j=0 \; \forall x_j\in \cnonmult_\div{J}^{\comp{\lm(\Sigma)}}(v)} \:= \: 
\varphi_{\beta_1,\ldots,\beta_n}(x_{i_1},\ldots,x_{i_{k_v}}),
\]
where $\{x_{i_1},\ldots,x_{i_{k_v}}\} = \cmult^{\comp{\lm(\Sigma)}}_{\div{J}}(v)$.
\end{theorem}  

\subsection{An algorithmic approach to solvability for monomial PDE systems}
With Theorem~\ref{Theorem:BoundaryConditions}, Janet gave a solution to the \emph{Cauchy problem} for a monomial PDE system $(\Sigma)$. To summarize Janet's approach, the Cauchy problem for the system $(\Sigma)$ can be solved by the following steps.
\begin{enumerate}[{\bf i)}]
\item If the set $\lm(\Sigma)$ of leading monomials of $(\Sigma)$ is complete,
\begin{itemize}
\item if all compatibility conditions in $(C_\Sigma)$ are satisfied, then the Cauchy problem admits a solution,
\item in the others cases, the system $(\Sigma)$ is incompatible. 
\end{itemize}
\item If the set $\lm(\Sigma)$ is not complete, then apply the step {\bf i)} to the completion of $\lm(\Sigma)$.
\end{enumerate}

Without the completeness property, a monomial PDE system $(\Sigma)$ may have infinitely many compatibility conditions. With the algorithmic approach introduced by Janet, these are reduced to a finite number of compatibility conditions of the form~\ref{E:ConditionI}. Indeed, it suffices to verify the conditions on a finite set that involutively generates the set $\lm(\Sigma)$ of leading monomials of the PDE system $(\Sigma)$.

\section{Janet's monomial order on derivatives}
\label{S:JanetMonomialOrderOnDerivative}

The main novelty in Janet's monograph \emph{Leçons sur les systèmes d'équations aux dérivées partielles}, \cite{Janet29}, published in 1929, is his treatment of the solvability problem of linear PDE systems defined by polynomial equations.
With the notion of order defined with principal and parametric derivative, he gave an algebraic characterization of complete integrability conditions of such systems. He also used this order to define a procedure that decides whether a given finite linear polynomial PDE system can be transformed into a completely integrable linear polynomial PDE system.
The solvability result presented in the previous section is based on a formulation of initial conditions in terms of complementary monomials. In this way, the partition~(\ref{E:MacaulayMonomial}) is essential in this approach. With a view to extending these construction to polynomial PDE systems, Janet considered an order on derivatives defined using the notions of principal and parametric derivative that take the partition~(\ref{E:MacaulayMonomial}) precisely into account.

\subsection{Principal and parametric derivatives}
In the 1929 monograph,~\cite{Janet29}, Janet extended  Theorem~\ref{Theorem:BoundaryConditions} on the Cauchy problem for monomial PDE systems to polynomial PDE systems. He considered PDE systems in analytic categories, namely those in which all unknown functions, coefficients and initial conditions are supposed to be analytic.
The analyticity hypothesis considered by Janet corresponds to the classical notion, namely a function is analytic on a neighborhood of a point if it admits an analytic expression as a convergent series on this neighborhood. 

Janet obtained a generalization of the Cauchy-Kowalevsky theorem by defining an order on the set of derivatives that is compatible with products. Orders with the property of respecting the products corresponds to the notion of \emph{monomial order}. Such an order was first used by Gauss  in the proof of the fundamental theorem of symmetric polynomials with the lexicographic order. Monomial orders appeared also in Paul Gordan's  proof of the Hilbert's basis theorem  published in \cite{Gordan1893}. Finally the notion of ideal with respect to lexicographic order appeared in the work of F. S. Macaulay in \cite{Macaulay27}.

As explained in Section~\ref{S:EmergenceFormalMethodsPDE}, the notion of principal and parametric derivative emerged in the works of Ch. Méray and Ch. Riquier in their work on solvability of linear PDE systems in the period 1890-1910. These notions were reformulated in an appropriate algebraic language by Janet. He presented a notion of order on derivatives in two steps. First, he considered a lexicographic order on derivatives already defined by \'E. Delassus, \cite{Delassus1896}, using the terminology of \emph{anteriority} and \emph{posteriority}. He wrote in~{\cite[pp. 308-309]{Janet21b}}:

\begin{quote}
\emph{
Convenons de dire que si deux dérivées $D$, $D'$ de même ordre ont pour indice respectivement $\alpha_1, \alpha_2, ..., \alpha_n$ ; $\alpha_1',\alpha_2',\ldots,\alpha_n'$, $D$ est postérieur ou antérieur à $D'$ suivant que la première des différences $\alpha_1 - \alpha_1', \alpha_2 -\alpha_2', \ldots, \alpha_n - \alpha_n'$ qui n'est pas nulle est positive ou négative. 
}
\end{quote}

\begin{quote}
\emph{Let us say that, if two derivatives $D, D'$ of the same order have the indices $\alpha_1, \alpha_2, ..., \alpha_n$ ; $\alpha_1',\alpha_2',\ldots,\alpha_n'$, respectively, $D$ is posterior or prior to $D'$ according as the first difference $\alpha_1 - \alpha_1', \alpha_2 -\alpha_2', \ldots, \alpha_n - \alpha_n'$ which is not zero is positive or negative.  
}
\end{quote}

\noindent Note that, Janet  reversed the definition of the notion of posteriority and anteriority.
Second, he defined the notion of \emph{principal derivative} and \emph{parametric derivative} with respect to the lexicographic order previously defined. He wrote in~\cite[pp. 312]{Janet21b}:

\begin{quote}
\emph{
Considérons, pour simplifier un peu l'exposition, un système à une seule fonction inconnue~$z$~; convenons que si $D$, $D'$ sont deux dérivées d'ordres différents $p$, $p'$, $D$ est postérieure ou antérieure à $D'$ suivant que $p$ est supérieur ou inférieur à $p'$ ; adoptons d'autre part pour les dérivées d'un même ordre le classement même qui a été défini plus haut. Soit $(E)$ l'une quelconque des équations que l'on peut déduire du système par dérivations et combinaisons~; résolvons-la par rapport à la dernière des dérivées qui y entrent effectivement~; ce mode de résolution distingue un certain nombre de dérivées de $z$, celles qui figurent dans les premiers membres : nous les appellerons principales, toutes les autres seront appelées paramétriques.
}
\end{quote}

\begin{quote}
\emph{Let us consider, to simplify the explanation a little, a system with only one unknown function $z$ ; say that, if $D$, $D'$ are two different derivatives of different order $p$, $p'$, $D$ is posterior or prior to $D'$ according as $p$ is greater than or less than $p'$ ;  and for the derivatives of the same order we adopt the same order as defined above. Let $(E)$ be any one of the equations that we can deduce from the system by derivations and combinations; we solve it with respect to the last derivatives contained in the equation ; this way of resolution distinguishes a certain number of derivatives of $z$, those which appear in the left hand side: we call them principal and all of the others are called parametric. 
}
\end{quote}

\subsection{Weighted parametric and principal derivatives}
\label{SS:ParametricPrincipalJanet}

The analysis of linear PDE systems is made with respect to a given order on the set of monomials associated to derivatives. In order to specify the order to the problem  being studied, Janet generalized the order defined using the previous notion of posteriority on derivatives by introducing some weights attached to the indeterminates of the system. This weighted order is inspired by the notion of \emph{cote} introduced by Ch. Riquier in \cite{Riquier93} and \'{E}. Delassus in \cite{Delassus1896}, as mentioned in the historical context Section~\ref{S:HistoricalContext}.
In his monograph, Janet first considered the degree lexicographic order, {\cite[\textsection 22]{Janet29}}, formulated as follows:
\begin{enumerate}[{\bf i)}]
\item for $\vert\alpha\vert\neq \vert\beta\vert$, the derivative $D^\alpha\varphi$ is called \emph{posterior} (resp. \emph{anterior}) to~$D^\beta\varphi$, if $\vert\alpha\vert>\vert\beta\vert$ (resp. $\vert\alpha\vert<\vert\beta\vert$),
\item for $\vert\alpha\vert=\vert\beta\vert$, the derivative $D^\alpha\varphi$ is called \emph{posterior} (resp. \emph{anterior}) to $D^\beta\varphi$ if the first non-zero difference 
\[ 
\alpha_n-\beta_n\; , \quad \alpha_{n-1}-\beta_{n-1}\; , \quad \ldots\quad, \; \alpha_1-\beta_1,
\]
is positive (resp. negative).
\end{enumerate}

Let us consider the following equation:
\begin{equation}
\label{EDP-2.2}
D\varphi\: = \:\sum_{i \in I}a_{i}D_i\varphi+f, 
\end{equation}
where $D$ and the $D_i$ are differential operators such that $D_i\varphi$ is anterior to $D\varphi$ for all $i$ in $I$. 
The derivative~$D\varphi$ and all its derivatives are called \emph{principal derivatives of Equation (\ref{EDP-2.2})}. All the other derivative of~$u$ are called \emph{parametric derivatives of Equation (\ref{EDP-2.2})}.

Further generalization of these order relations were given by Janet by formulating a new notion of \emph{cote}, that corresponds to a parametrization of a weight order defined as follows.
Let us fix a positive integer $s$. We define a matrix of \emph{weight}
\[
C \: = \: \left[
\begin{tabular}{cccccc} 
$C_{1,1}$  & $\ldots$ &  $C_{n,1}$ \\
$\vdots$ & &$\vdots$  \\
$C_{1,s}$ &  $\ldots$ & $C_{n,s}$ \\
\end{tabular}
\right]
\]
that associates to each variable~$x_i$ non negative integers $C_{i,1}, \ldots, C_{i,s}$, called the \emph{$s$-weights} of~$x_i$. This notion was called \emph{cote} by Janet in {\cite[\textsection 22]{Janet29}} following the terminology introduced by Ch. Riquier,~\cite{Riquier10}. 
For each derivative $D^\alpha\varphi$, with $\alpha=(\alpha_1,\ldots,\alpha_n)$ of an analytic function $\varphi$, we associate a \emph{$s$-weight}~$\Gamma(C)=(\Gamma_1,\ldots,\Gamma_s)$ where the $\Gamma_k$ are defined by 
\[ 
\Gamma_k \: = \:\sum_{i=1}^n \alpha_i C_{i,k}.
\]
Given two monomial partial differential operators $D^\alpha$ and $D^\beta$, we say that $D^\alpha\varphi$ is \emph{posterior} (resp. \emph{anterior}) to $D^\beta\varphi$ with respect to a weigh matrix $C$ if 
\begin{enumerate}[{\bf i)}]
\item $\vert\alpha\vert\neq \vert\beta\vert$ and $\vert\alpha\vert>\vert\beta\vert$ (resp. $\vert\alpha\vert<\vert\beta\vert$),
\item otherwise $\vert\alpha\vert=\vert\beta\vert$ and the first non-zero difference  
\[ 
\Gamma_1-\Gamma_1',  \quad \Gamma_2-\Gamma_2'\; , \quad \ldots\quad, \;  \Gamma_s-\Gamma_s',
\]
is positive (resp. negative). 
\end{enumerate}
In this way, we define an order on the set of monomial partial derivatives, called \emph{weight order}.
Note that, this notion generalizes the above lexicographic order defined by Janet, that corresponds to the case~$C_{i,k}=\delta_{i+k,n+1}$.

\subsection{Complete higher-order finite linear PDE systems}
\label{SS:CompleteHigherOrderFiniteLineraPDESystems}

In {\cite[\textsection 39]{Janet29}}, Janet studied the solvability of the following PDE system of one unknown function~$\varphi$ in which each equation is of the following form:\begin{equation}
\label{Equation:HigherOrderPDE}
(\Sigma)\qquad
D_i\varphi\: = \:\sum_{j}a_{i,j}D_{i,j}\varphi, \quad i \in I,
\end{equation}
where all the functions $a_{i,j}$ are supposed analytic in a neighborhood of a point $P$  in $\mathbb{C}^n$, and each equation is supposed to satisfy the following two conditions:
\begin{enumerate}[{\bf i)}]
\item $D_{i,j}\varphi$ is anterior to $D_i\varphi$, for any $i$ in $I$,
\item all the $D_i$'s for $i$ in $I$ are distinct.
\end{enumerate}
He defined the notion of principal derivative for such a system by setting: the derivatives $D_i\varphi$, for $i$ in $I$, and all their derivatives, are called \emph{principal derivatives} of the PDE system $(\Sigma)$ given in (\ref{Equation:HigherOrderPDE}). Any other derivative of~$\varphi$ is called \emph{parametric derivative}. 
In this way, to the set of operators $D_i$ for $i$ in $I$,  he associated a set $\lm(\Sigma)$ of monomials through the morphism $\Phi$ defined Section~\ref{SS:MonomialPartialDifferentialEquationsSystems}.
The PDE system $(\Sigma)$ is then said to be \emph{complete} if the set of monomials~$\lm(\Sigma)$ is complete. Note that in \cite{Janet20}, Janet introduced a completion procedure that transforms a finite linear PDE system into an equivalent complete linear PDE system.

By definition, the set of principal derivatives corresponds to the multiplicative cone of $\lm(\Sigma)$. Hence, when the system $(\Sigma)$ is complete, the set of principal derivatives corresponds to the involutive cone of $\lm(\Sigma)$. Having the partition 
\[
\Mr(x_1,\ldots,x_n) \: = \: \cone(\lm(\Sigma)) \amalg {\cone}^{\complement}_{\div{J}}(\comp{\lm(\Sigma)}),
\]
the set of parametric derivatives of the complete system $(\Sigma)$ corresponds to the involutive cone of the set $\comp{\lm(\Sigma)}$ of complementary monomials of $\lm(\Sigma)$.
To a monomial $x^{\beta}$ in~$\comp{\lm(\Sigma)}$, with $\beta=(\beta_1,\ldots,\beta_n)$ in $\Nb^n$ and 
\[
\cmult_{\div{J}}^{\comp{\lm(\Sigma)}}(x^{\beta}) \: = \: \{x_{i_1},\ldots,x_{i_{k_\beta}}\}, 
\]
we associate an arbitrary analytic function $\varphi_{\beta}(x_{i_1},\ldots,x_{i_{k_\beta}})$.
Using these functions, Janet defined a \emph{initial condition}:
\[
(C_\beta) \qquad \left. D^\beta \varphi \right|_{x_j=0 \; \forall x_j\in \cnonmult_{\div{J}}^{\comp{\lm(\Sigma)}}(x^{\beta})}\: = \: 
\varphi_{\beta}(x_{i_1},\ldots,x_{i_{k_\beta}}).
\]

\begin{theorem}[{\cite[\textsection 39]{Janet29}}]
\label{thm_exist-EDP1}
If the PDE system $(\Sigma)$ in~(\ref{Equation:HigherOrderPDE}) is complete, then it admits at most one analytic solution satisfying the initial condition 
\begin{equation}
\label{Equation:BoundaryCondition}
\{\, (C_\beta) \;|\; x^\beta \in \comp{\lm(\Sigma)}\,\}.
\end{equation}
\end{theorem}

Note that this result does not prove the existence of a solution of the PDE system $(\Sigma)$. The existence of solutions will be discussed in Section~\ref{SS:CompletelyIntegrableSystems}.

As we observed, the values of the parametric derivatives completely determine the initial condition~(\ref{Equation:BoundaryCondition})
That is, these derivatives parameterize the space of solutions of the differential equation~(\ref{Equation:HigherOrderPDE}). This observation suggests the origin of the terminology \emph{parametric derivative} introduced by Ch. Méray and Ch. Riquier. 

\subsection{Linear PDE systems for several unknown functions}
Janet extended the construction of initial conditions given above for one unknown function to linear PDE systems on $\mathbb{C}^n$ with several unknown functions using a weight order. Consider a linear PDE system of $m$ unknown analytic functions~$\varphi^1,\ldots, \varphi^m$ of the following form
\begin{equation}
\label{Equation:CanonicalSystem0}
(\Sigma) \qquad D^\alpha \varphi^r\: = \:\sum_{\substack{(\beta,s) \in \mathbb{N}^n\times\{1,2,\ldots,m\}}}
a_{\alpha,\beta}^{r,s}D^\beta \varphi^s, \qquad \alpha \in I^r,
\end{equation}
for $1\leq r\leq m$, where $I^r$ is a finite subset of $\mathbb{N}^n$ and the $a_{\alpha,\beta}^{r,s}$ are analytic functions. He defined a weight order in such a way that the system~(\ref{Equation:CanonicalSystem0}) can be expressed in the form 
\begin{equation}
\label{Equation:CanonicalSystem1}
(\Sigma) \qquad D^\alpha \varphi^r\: = \:\sum_{\substack{(\beta,s) \in \mathbb{N}^n\times\{1,2,\ldots,m\} \\ D^\beta \varphi^s \wostrict D^\alpha \varphi^r}}
a_{\alpha,\beta}^{r,s}D^\beta \varphi^s, \qquad \alpha \in I^r,
\end{equation}
allowing him to formulate the notion of completeness of the system $(\Sigma)$. Let $\lm_{\wo}(\Sigma,\varphi^r)$ be the set of monomials associated to leading derivatives $D^\alpha$ of all PDE in $(\Sigma)$ such that $\alpha$ belongs to $I^r$. The PDE system $(\Sigma)$ is called \emph{complete} with respect to $\wo$, if for any $1\leq r \leq m$, $\lm_{\wo}(\Sigma,\varphi^r)$ is complete as a set of monomials.

The question is to determine under which conditions the system $(\Sigma)$  in~(\ref{Equation:CanonicalSystem1}) admits a solution for any given initial condition. 
We suppose that $(\Sigma)$ is complete, hence the set of monomials  $\lm_{\wo}(\Sigma,\varphi^r)=\{x^\alpha \;|\; \alpha \in I^r\}$, which we will denote by $\Ur_r$, is complete for all $1\leq r\leq m$. 
The initial conditions for which the system admits at most one solution are parametrized by the set~$\comp{\Ur_r}$ of complementary monomials of the set of monomials $\Ur_r$. Explicitly, for $1\leq r \leq m$, to a monomial $x^{\beta}$ in $\comp{\Ur_r}$, with $\beta$ in $\Nb^n$ and $\cmult_{\div{J}}^{\comp{\Ur_r}}(x^{\beta}) = \{x_{i_1},\ldots,x_{i_{k_r}}\}$, we associate an arbitrary analytic function 
\[
\varphi_{\beta,r}(x_{i_1},\ldots,x_{i_{k_r}}).
\]
Formulating \emph{initial condition} as the following data:
\[
(C_{\beta,r}) \qquad \left. D^{\beta} \varphi^r \right|_{x_j=x_j^0 \; \forall x_j\in \cnonmult_{\div{J}}^{\comp{\Ur_r}}(x^{\beta})}\: = \: 
\varphi_{\beta,r}(x_{i_1},\ldots,x_{i_{k_r}}),
\]
we set the \emph{initial condition} of the system $(\Sigma)$ in~(\ref{Equation:CanonicalSystem0}) to be the following set
\begin{equation}
\label{Equation:BoundaryCondition2}
\underset{1\leq r \leq m}{\bigcup}\{\, C_{\beta,r} \;|\; x^{\beta} \in \comp{\Ur_r}\,\}.
\end{equation}
Explaining that the proof is similar to the proof of Theorem~\ref{thm_exist-EDP1}, 
Janet announced the following result.

\begin{theorem}[{\cite[\textsection 40]{Janet29}}]
\label{thm_exist-EDP2}
If the PDE system $(\Sigma)$ in~(\ref{Equation:CanonicalSystem1}) is complete with respect to a weight order~$\wo$, then it admits at most one analytic solution satisfying the initial condition (\ref{Equation:BoundaryCondition2}).
\end{theorem}

\subsection{Completely integrable systems}
\label{SS:CompletelyIntegrableSystems}

Given $1\leq r\leq m$ and  $\alpha \in I^r$, let $x_i$ be in $\nonmult_{\div{J}}^{\Ur_r}(x^\alpha)$ a non-multiplicative variable. Let us differentiate the equation
\[ 
D^\alpha \varphi^r\: = \:\sum_{\substack{(\beta,s) \in \mathbb{N}^n\times\{1,2,\ldots,m\} \\ D^\beta \varphi^s \wostrict  D^\alpha \varphi^r}}
a_{\alpha,\beta}^{r,s}D^\beta \varphi^s
\]
by the partial derivative $\Phi(x_i)=\frac{\partial}{\partial x_i}$. We obtain the following PDE
\begin{equation}
\label{Equation:CanonicalSystem}
\Phi(x_i)(D^\alpha \varphi^r)\: = \:
\sum_{\substack{(\beta,s) \in \mathbb{N}^n\times\{1,2,\ldots,m\} \\ D^\beta \varphi^s \wostrict D^\alpha \varphi^r}} \left(\frac{\partial a_{\alpha,\beta}^{r,s}}{\partial x_i}D^\beta \varphi^s+a_{\alpha,\beta}^{r,s}\Phi(x_i)(D^\beta \varphi^s)\right). 
\end{equation}
Using the system $(\cone_{\div{J},\wo}(\Sigma))$, we can rewrite the PDE~(\ref{Equation:CanonicalSystem}) into a PDE formulated in terms of parametric derivatives and independent variables. The set of monomials $\Ur_r$ being complete, there exists $\alpha'$  in $\Nb^n$ with $x^{\alpha'}$ in $\Ur_r$ and $u$  in~$\mathcal{M}(\mult_{\div{J}}^{\Ur_r}(x^{\alpha'}))$ such that $x_ix^{\alpha}=ux^{\alpha'}$. Then $\Phi(x_i)D^\alpha=\Phi(u)D^{\alpha'}$, and as a consequence we obtain the following equation
\begin{equation}
\label{Equation:CanonicalSystem2}
\sum_{\substack{(\beta,s) \in \mathbb{N}^n\times\{1,2,\ldots,m\} \\ D^\beta \varphi^s \wostrict D^\alpha \varphi^r}}
\left(\frac{\partial a_{\alpha,\beta}^{r,s}}{\partial x_i}D^\beta \varphi^s+a_{\alpha,\beta}^{r,s}\Phi(x_i)(D^\beta \varphi^s)\right)
\: = \:
\sum_{\substack{(\beta',s) \in \mathbb{N}^n\times\{1,2,\ldots,m\} \\ D^{\beta'} \varphi^s \wostrict D^{\alpha'} \varphi^r}}
\Phi(u)(a_{\alpha',\beta'}^{r,s}D^{\beta'} \varphi^s).
\end{equation}
Using equations of the system $(\cone_{\div{J},\wo}(\Sigma))$, we replace all  principal derivatives in the equation (\ref{Equation:CanonicalSystem2}) by parametric derivatives and independent variables. The order $\wo$ being well-founded, this process will terminate. Moreover, when the PDE system $(\Sigma)$ is complete, this reduction process is confluent in the sense that any transformation of an equation~(\ref{Equation:CanonicalSystem2}) ends on a unique $\div{J}$-normal form.
This set of $\div{J}$-normal forms is denoted by $\integralcond_{\div{J},\wo}(\Sigma)$. 

The system $(\Sigma)$ being complete, any equation~(\ref{Equation:CanonicalSystem2}) is reduced to a unique normal form. Such a normal form allows us to judge whether a given integrability condition is trivial or not.
Recall that the parametric derivatives correspond to the initial conditions. Hence, a non-trivial relation in 
$\integralcond_{\div{J},\wo}(\Sigma)$ provides a non-trivial relation among the initial conditions. In this way, we can decide whether the system $(\Sigma)$ is completely integrable or not.
A complete linear PDE system $(\Sigma)$ of the form~(\ref{Equation:CanonicalSystem1}) is said to be \emph{completely integrable} if it admits an analytic solution for any given initial condition~(\ref{Equation:BoundaryCondition2}).

\begin{theorem}[{\cite[\textsection 42]{Janet29}}]
\label{Theorem:CaracterizationCompleteIntegrability}
Let $(\Sigma)$ be a complete finite linear PDE system of the form~(\ref{Equation:CanonicalSystem1}). Then the system $(\Sigma)$ is completely integrable if and only if any relation in $\integralcond_{\div{J},\wo}(\Sigma)$ is a trivial identity.
\end{theorem}

A proof of this result is given in {\cite[\textsection 43]{Janet29}}.
Note that the latter condition is equivalent to saying that any relation~(\ref{Equation:CanonicalSystem2}) is an algebraic consequence of a PDE equation of the system $(\cone_{\div{J},\wo}(\Sigma))$.

\subsection{The notion of canonical PDE system}
\label{SS:NotionCanonicalPDESystems}

In {\cite[\textsection 46]{Janet29}} Janet introduced the notion of canonical linear PDE system. A canonical system is a normal form with respect to a weight order on derivatives, and satisfying some analytic conditions, allowing an extension of the Cauchy-Kowalevsky theorem.
Janet gave a procedure which transforms a finite linear PDE system with several unknown functions into an equivalent linear PDE system that is either in canonical form or in an incompatible system. 
Janet formulated its procedure as follows, {\cite[\textsection 46]{Janet29}},

\begin{quote}
\emph{
Adoptons pour les variables indépendantes et les fonctions inconnues un système de cotes tel que chacune des classes qui en résultent ne contienne qu'un élément ;
}

\quad [$\ldots$]

\emph{
\'Etant donné un système quelconque donné $S$, comprenant un nombre fini d'équations, considérons la dernières $\Delta$, des dérivées qui y entrent, c'est-à-dire celle qui est postérieure à toutes les autres et résolvons par rapport à elle une des équations du système qui la contiennent ; portons l'expression trouvée dans les autres équations ; traitons le système obtenu qui ne contient pas $\Delta_1$ comme nous avons traité le système primitif, et ainsi de suite. Nous obtiendrons finalement un système $(\Sigma)$ d'équations résolues, chacune ne contenant dans son second membre que des dérivées antérieures à son premier membre, les premiers membres étant tous différents.}

\emph{
Formons les conditions d'intégrabilité complète $(C)$ du système obtenu. Nous obtiendrons des relations en nombre fini, ne contenant que les variables indépendantes et les dérivées paramétriques, qui, si le système n'est pas complètement intégrable, ne sont pas toutes des identités.
}

\emph{
Résolvons ces relations comme nous avons résolu celles du système primitivement donné $S$, et joignons les équations obtenues aux équations $(\Sigma)$. Nous obtenons un système $(\Sigma')$ formé encore d'équations résolues, chacune ne contenant dans son second membre que des dérivées antérieures à son premier membre, les premiers membres étant tous différents. Les premiers membres $(\Delta')$ de $(\Sigma')$ comprennent les premiers membres $(\Delta)$ de $(\Sigma)$ et des dérivées qui ne sont dérivées d'aucun des $(\Delta)$ puisque ce sont des dérivées paramétriques pour $(\Sigma)$. Nous traiterons $\Sigma'$ comme nous avons traité $\Sigma$, et ainsi de suite.
}

\emph{
Je dis que l'opération ne peut se répéter qu'un nombre fini de fois.
}

\end{quote}

\begin{quote}
\emph{For the independent variables and the unknown functions, adopt a system of $\ll$\,cote\,$\gg$ such that each class that is defined with respect to this system contains only one element ;
}

\quad [$\ldots$]

\emph{
Given any system $S$, containing a finitely number of equations, consider the last $\Delta$, the derivatives contained in equations, namely, the derivatives that is posterior to all other derivatives and solve one of the equations, containing a derivative, with respect to the derivative ; keep the expression found in the other equations, treat the obtained system which does not contain $\Delta_1$ as we have treated the primitive system, and so on. Finally, we obtain a system $(\Sigma)$ of solved equations, each equation
that contains in its right hand side only prior derivatives with respect to its left hand side, the terms of the left hand sides of all equations are different.
}

\emph{We form the complete integrability conditions $(C)$ of the obtained system.
We obtain a finite number of relations, which only contain the independent variables and the parametric derivatives, where, if the system is not completely integrable, not all of the relations are identities.  
}

\emph{Solve these relations as we solved for the primitively given system $S$, and join the obtained equations to the equations $(\Sigma)$. 
We will obtain a system $(\Sigma')$ formed by solved equations, where in the second member,  each equation containing only the prior derivatives to its first member and the first members are  all different. 
The first members $(\Delta')$ of $(\Sigma')$ contain the first members $(\Delta)$ of $(\Sigma)$ and the derivatives which are not derivatives of $(\Delta)$ because they are parametric derivatives for $(\Sigma)$. We will treat $\Sigma'$ as we treated $\Sigma$, and so on. 
}

\emph{
I claim that the operation can be repeated only finitely many times. 
}

\end{quote}

Let us formulate in the modern language explained in this article the notion of canonical form so obtained by Janet. Given a fixed weight order~$\wo$, we suppose that each equation of a finite linear PDE system $(\Sigma)$ can be expressed in the following form
\[
(\Sigma^{(\alpha,r)})
\qquad
D^\alpha \varphi^r\: = \:\sum_{\substack{(\beta,s) \in \mathbb{N}^n\times\{1,2,\ldots,m\} \\ D^\beta \varphi^s \wostrict D^\alpha \varphi^r}}
a_{(\beta,s)}^{(\alpha,r)}D^\beta \varphi^s.
\]
The support of the equation $(\Sigma^{(\alpha,r)})$ is defined by
\[
\mathrm{Supp}( \Sigma^{(\alpha,r)}) \: = \: \{\,(\beta,s)\;|\;a_{(\beta,s)}^{(\alpha,r)} \neq 0\;\}.
\]

For $1\leq r \leq m$, consider the set of monomials~$\lm_{\wo}(\Sigma,\varphi^r)$ corresponding to leading derivatives, that is monomials $x^{\alpha}$ such that $(\alpha,r)$ belongs to $I$.
The system $(\Sigma)$ is said to be 
\begin{enumerate}[{\bf i)}]
\item \emph{$\div{J}$-left-reduced with respect to $\wo$} if for any $(\alpha,r)$ in $I$ there is no $(\alpha',r)$ in $I$ and non-trivial monomial $x^\gamma$ in $\Mr(\mult_{\div{J}}^{\lm_{\wo}(\Sigma,\varphi^r)}(x^{\alpha'}))$ such that $x^\alpha=x^\gamma x^{\alpha'}$,
\item \emph{$\div{J}$-right-reduced with respect to $\wo$} if, for any $(\alpha,r)$ in $I$ and any $(\beta,s)$ in $\mathrm{Supp}( \Sigma^{(\alpha,r)})$, there is no~$(\alpha',s)$ in 
$I$ and non-trivial monomial $x^\gamma$ in $\Mr(\mult_{\div{J}}^{\lm_{\wo}(\Sigma,\varphi^r)}(x^{\alpha'}))$ such that $x^\beta=x^\gamma x^{\alpha'}$,
\item \emph{$\div{J}$-autoreduced with respect to $\wo$} if it is both $\div{J}$-left-reduced and $\div{J}$-right-reduced with respect to $\wo$.
\end{enumerate}

A PDE system $(\Sigma)$ is said to be \emph{$\div{J}$-canonical with respect a weight order~$\wo$} if it satisfies the following five conditions
\begin{enumerate}[{\bf i)}]
\item it consists of finitely many equations and each equation can be expressed in the following form
\[
D^\alpha \varphi^r\: = \:\sum_{\substack{(\beta,s) \in \mathbb{N}^n\times\{1,2,\ldots,m\} \\ D^\beta \varphi^s \wostrict D^\alpha \varphi^r}}
a_{(\beta,s)}^{(\alpha,r)}D^\beta \varphi^s,
\]
\item the system $(\Sigma)$ is $\div{J}$-autoreduced with respect to $\wo$,
\item the system $(\Sigma)$ is complete,
\item the system $(\Sigma)$ is completely integrable,
\item the coefficients $a_{(\beta,s)}^{(\alpha,r)}$ of the equations in {\bf i)} and the initial conditions of $(\Sigma)$ are analytic.
\end{enumerate}
Under these assumptions, the system $(\Sigma)$ admits a unique analytic solution satisfying appropriate initial conditions parametrized by complementary monomials.
In his monograph~\cite{Janet29}, Janet did not mention the correctness of the procedures that he introduced in order to reduce a finite linear PDE system to a canonical form. We refer the reader to \cite{IoharaMalbos19ACM} for a more complete account on the Janet procedure.

\begin{small}
\bibliographystyle{plain}
\bibliography{biblioCURRENT}

\def\cprime{$'$}
\begin{thebibliography}{10}

\bibitem{AubinGispertGoldstein14}
D.~Aubin, H.~Gispert, and C.~Goldstein.
\newblock The total war of {P}aris mathematicians.
\newblock In {\em The war of guns and mathematics}, volume~42 of {\em Hist.
  Math.}, pages 125--177. Amer. Math. Soc., Providence, RI, 2014.

\bibitem{AubinGoldstein14}
D.~Aubin and C.~Goldstein.
\newblock Placing {W}orld {W}ar {I} in the history of mathematics.
\newblock In {\em The war of guns and mathematics}, volume~42 of {\em Hist.
  Math.}, pages 1--55. Amer. Math. Soc., Providence, RI, 2014.

\bibitem{Bergman78}
G.~M. Bergman.
\newblock The diamond lemma for ring theory.
\newblock {\em Adv. in Math.}, 29(2):178--218, 1978.

\bibitem{Bourlet91}
C.~Bourlet.
\newblock Sur les \'equations aux d\'eriv\'ees partielles simultan\'ees qui
  contiennent plusieurs fonctions inconnues.
\newblock {\em Annales scientifiques de l'\'Ecole Normale Sup\'erieure}, 3e
  s{\'e}rie, 8:3--63, 1891.

\bibitem{BC3G91}
R.~L. Bryant, S.~S. Chern, R.~B. Gardner, H.~L. Goldschmidt, and P.~A.
  Griffiths.
\newblock {\em Exterior differential systems}, volume~18 of {\em Mathematical
  Sciences Research Institute Publications}.
\newblock Springer-Verlag, New York, 1991.

\bibitem{Buchberger65}
B.~Buchberger.
\newblock {\em {Ein Algorithmus zum Auffinden der Basiselemente des
  Restklassenringes nach einem nulldimensionalen Polynomideal (An Algorithm for
  Finding the Basis Elements in the Residue Class Ring Modulo a Zero
  Dimensional Polynomial Ideal)}}.
\newblock PhD thesis, Mathematical Institute, University of Innsbruck, Austria,
  1965.
\newblock English translation in J. of Symbolic Computation, Special Issue on
  Logic, Mathematics, and Computer Science: Interactions. Vol. 41, Number 3-4,
  Pages 475--511, 2006.

\bibitem{BuraliForti1897}
C.~Burali-Forti.
\newblock {\em Introduction \`{a} la {G}\'{e}om\'{e}trie {D}iff\'{e}rentielle
  suivant la {M}\'{e}thode de {H}. {G}rassmann}.
\newblock Gauthier-Villars, 1897.

\bibitem{Cartan1899}
{\'E}.~{Cartan}.
\newblock {Sur certaines expressions diff\'erentielles et le probl\`eme de
  Pfaff.}
\newblock {\em {Ann. Sci. \'Ec. Norm. Sup\'er. (3)}}, 16:239--332, 1899.

\bibitem{Cogliati11}
A.~Cogliati.
\newblock On the genesis of the {C}artan-{K}\"ahler theory.
\newblock {\em Arch. Hist. Exact Sci.}, 65(4):397--435, 2011.

\bibitem{Darboux82}
G.~{Darboux}.
\newblock {Sur le probl\`eme de Pfaff.}
\newblock {\em {Bull. Sci. Math. et Astro.}}, 6:14--36, 49--68, 1882.

\bibitem{Dedekind77}
R.~{Dedekind}.
\newblock {Sur la th\'eorie des nombres entiers alg\'ebriques.}
\newblock {\em {Darboux Bull. (2)}}, 1:17--41, 69--92, 144--164, 207--248,
  1877.

\bibitem{Delassus1896}
\'{E}. Delassus.
\newblock Extension du th\'eor\`eme de {C}auchy aux syst\`emes les plus
  g\'en\'eraux d'\'equations aux d\'eriv\'ees partielles.
\newblock {\em Ann. Sci. \'Ecole Norm. Sup. (3)}, 13:421--467, 1896.

\bibitem{Dickson13}
L.~E. Dickson.
\newblock Finiteness of the {O}dd {P}erfect and {P}rimitive {A}bundant
  {N}umbers with {$n$} {D}istinct {P}rime {F}actors.
\newblock {\em Amer. J. Math.}, 35(4):413--422, 1913.

\bibitem{Dugas50}
R.~Dugas.
\newblock {\em Histoire de la m\'ecanique}.
\newblock Editions du Griffon, Neuchatel, 1950.

\bibitem{Forsyth90}
A.~R. {Forsyth}.
\newblock {Theory of differential equations. Part I. Exact equations and
  Pfaff's problem.}
\newblock {Cambridge. University Press. XIII + 340 S. $8^{\circ}$}, 1890.

\bibitem{Frobenius77}
G.~Frobenius.
\newblock Ueber das {P}faffsche {P}roblem.
\newblock {\em J. Reine Angew. Math.}, 82:230--315, 1877.

\bibitem{Gerdt97}
V.~P. Gerdt.
\newblock Gr\"obner bases and involutive methods for algebraic and differential
  equations.
\newblock {\em Math. Comput. Modelling}, 25(8-9):75--90, 1997.
\newblock Algorithms and software for symbolic analysis of nonlinear systems.

\bibitem{Gordan1893}
P.~{Gordan}.
\newblock {Ueber einen Satz von Hilbert.}
\newblock {\em {Math. Ann.}}, 42:132--142, 1893.

\bibitem{Goursat22}
\'E. {Goursat}.
\newblock {Le\c cons sur le probl\`eme de Pfaff.}
\newblock {Paris: J. Hermann. VIII u. 386 S. $8^\circ$}, 1922.

\bibitem{Grassmann44}
H.~E. {Grassmann}.
\newblock {Die lineale Ausdehnungslehre}.
\newblock {Leipzig: Verlag von Otto Wigand. 324 p. }, 1844.

\bibitem{Griffiths83}
P.~A. Griffiths.
\newblock {\em Exterior differential systems and the calculus of variations},
  volume~25 of {\em Progress in Mathematics}.
\newblock Birkh\"auser, Boston, Mass., 1983.

\bibitem{Grobner37}
W.~Gr{\"o}bner.
\newblock \"{U}ber das macaulaysche inverse system und dessen bedeutung f\"ur
  die theorie der linearen differentialgleichungen mit konstanten
  koeffizienten.
\newblock {\em Abh. Math. Sem. Univ. Hamburg}, 12(1):127--132, 1937.

\bibitem{Gunther13b}
N.~G{\"u}nther.
\newblock {\"Uber die kanonische Form der Systeme kanonischer homogener
  Gleichungen.}
\newblock {Samml. des Inst. der Verkehrswege 82, 22 S. l~{\ss}}, 1913.

\bibitem{Gunther41}
N.~M. {Gunther}.
\newblock {Sur les modules des formes alg\'ebriques.}
\newblock {Trudy Tbilis. Mat. Inst. 9, 97-206}, 1941.

\bibitem{Hawkins05}
T.~Hawkins.
\newblock Frobenius, {C}artan, and the problem of {P}faff.
\newblock {\em Arch. Hist. Exact Sci.}, 59(4):381--436, 2005.

\bibitem{Hilbert1890}
D.~Hilbert.
\newblock Ueber die {T}heorie der algebraischen {F}ormen.
\newblock {\em Math. Ann.}, 36(4):473--534, 1890.

\bibitem{Hironaka64}
H.~Hironaka.
\newblock Resolution of singularities of an algebraic variety over a field of
  characteristic zero. {I}, {II}.
\newblock {\em Ann. of Math. (2) 79 (1964), 109--203; ibid. (2)}, 79:205--326,
  1964.

\bibitem{IoharaMalbos19ACM}
Kenji {Iohara} and Philippe {Malbos}.
\newblock {From analytical mechanics problems to rewriting theory through {M}.
  {J}anet's work.}
\newblock In {\em {Two algebraic byways from differential equations: Gr\"obner
  bases and quivers}}, pages 3--74. Algorithms and Computation in Mathematics
  28, Springer, 2020.

\bibitem{Jacobi27b}
C.G.J. {Jacobi}.
\newblock {Ueber die Pfaffsche Methode, eine gew\"{o}nliche line\"{a}re
  {D}ifferentialgleichung zwischen $2n$ {V}ariablen durch ein {S}ystem von $n$
  {G}leichungen zu integriren}.
\newblock {\em {J. Reine Angew. Math.}}, 2:347--357, 1827.

\bibitem{Janet13b}
M.~{Janet}.
\newblock {Existence et d\'etermination univoque des solutions des syst\`emes
  d'\'equations aux deriv\'ees partielles.}
\newblock {\em {C. R. Acad. Sci., Paris}}, 157:697--700, 1913.

\bibitem{Janet13a}
M.~{Janet}.
\newblock {Sur les caract\'eristiques des syst\`emes d'\'equations aux
  deriv\'ees partielles.}
\newblock {\em {C. R. Acad. Sci., Paris}}, 156:118--121, 1913.

\bibitem{Janet20PhD}
M.~{Janet}.
\newblock {\em Sur les syst{\`e}mes d'{\'e}quations aux d{\'e}riv{\'e}es
  partielles}.
\newblock PhD thesis, Faculté des sciences de Paris, 6 1920.
\newblock Gauthier-Villars, Paris.

\bibitem{Janet20}
M.~{Janet}.
\newblock Sur les syst{\`e}mes d'{\'e}quations aux d{\'e}riv{\'e}es partielles.
\newblock {\em Journal de mathématiques pures et appliquées}, 8(3):65--151,
  1920.

\bibitem{Janet20a}
M.~{Janet}.
\newblock {Sur les syst\`emes d'\'equations aux d\'eriv\'ees partielles.}
\newblock {\em {C. R. Acad. Sci., Paris}}, 170:1101--1103, 1920.

\bibitem{Janet21a}
M.~{Janet}.
\newblock {Sur la recherche g\'en\'erale des fonctions primitives \`a $n$
  variables.}
\newblock {\em {Bull. Sci. Math., II. S\'er.}}, 45:238--248, 1921.

\bibitem{Janet22a}
M.~{Janet}.
\newblock {Les caract\`eres des modules de formes et les syst\`emes
  d'\'equations aux d\'eriv\'ees partielles.}
\newblock {\em {C. R. Acad. Sci., Paris}}, 174:432--434, 1922.

\bibitem{Janet24}
M.~{Janet}.
\newblock {Les modules de formes alg\'ebriques et la th\'eorie g\'en\'erale des
  syst\`emes diff\'erentielles.}
\newblock {\em {Ann. Sci. \'Ec. Norm. Sup\'er. (3)}}, 41:27--65, 1924.

\bibitem{Janet21b}
M.~{Janet}.
\newblock {Les travaux r\'ecents sur le degr\'e d'ind\'etermination des
  solutions d'un syst\`eme diff\'erentiel.}
\newblock {\em {Bull. Sci. Math., II. S\'er.}}, 49:307--320, 332--344, 1925.

\bibitem{Janet-ICM24}
M.~{Janet}.
\newblock {Sur les syst\`emes lin\'eaires d'hypersurfaces.}
\newblock {Proceedings Congress Toronto 1, 835-841}, 1928.

\bibitem{Janet29}
M.~{Janet}.
\newblock {Le\c{c}ons sur les syst\`emes d'\'equations aux d\'eriv\'ees
  partielles.}
\newblock {VIII + 124 p. Paris, Gauthier-Villars (Cahiers scientifiques
  publi\'es sous la direction de {\it G.~Julia}, fasc.~IV.)}, 1929.

\bibitem{Janet-ICM32}
M.~{Janet}.
\newblock {D\'etermination explicite de certains minima .}
\newblock {Verhandlungen Kongre{\ss } Z\"urich 1932, 2, 111-113}, 1932.

\bibitem{Janet-ICM36}
M.~{Janet}.
\newblock {Sur les syst\`emes de deux \'equations aux d\'eriv\'ees partielles
  \`a deux fonctions inconnues.}
\newblock {C. R. Congr. internat. Math., Oslo 1936, 2, 61-62}, 1936.

\bibitem{Jordan1870}
C.~{Jordan}.
\newblock {\em {Trait\'e des substitutions et des \'equations alg\'ebriques.}}
\newblock Paris: \'Editions J. Gabay, r\'eimpression du orig. 1870 edition,
  1989.

\bibitem{Kahler34}
E.~{K{\"a}hler}.
\newblock {Einf\"uhrung in die Theorie der Systeme von
  Differentialgleichungen.}
\newblock {(Hamburg. Math. Einzelschr. 16) Leipzig, Berlin: B. G. Teubner IV,
  80 S}, 1934.

\bibitem{Katz85}
V.-J. Katz.
\newblock Differential forms---{C}artan to de {R}ham.
\newblock {\em Arch. Hist. Exact Sci.}, 33(4):321--336, 1985.

\bibitem{Koenig1903}
J.~{K\"onig}.
\newblock {Einleitung in die allgemeine Theorie der algebraischen
  Gr\"o{\ss}en.}
\newblock {Leipzig: B. G. Teubner. X u. 552 S. \(8^{\circ}\) (1903).}, 1903.

\bibitem{Kowalevsky75}
S.~{Kowalevsky}.
\newblock {Zur Theorie der partiellen Differentialgleichungen.}
\newblock {\em {J. Reine Angew. Math.}}, 80:1--32, 1875.

\bibitem{Kronecker82}
L.~{Kronecker}.
\newblock {Grundz\"uge einer arithmetischen Theorie der algebraischen
  Gr\"ossen. (Festschrift zu Herrn Ernst Eduard Kummers f\"unfzigj\"ahrigem
  Doctor-Jubil\"aum, 10 September 1881).}
\newblock {\em {J. Reine Angew. Math.}}, 92:1--122, 1882.

\bibitem{Lagrange72}
J.-L. {Lagrange}.
\newblock Sur l'int\'{e}gration des \'{e}quations \`{a} diff\'{e}rences
  partielles du premier ordre.
\newblock {\em M\'{e}m. Acad. Sci. et Belles-Lettres de Berlin}, pages
  353--372, 1772.

\bibitem{Lagrange88}
J.-L. {Lagrange}.
\newblock {\em Méchanique Analitique}.
\newblock Desaint, 1788.

\bibitem{LejeuneDirichlet94}
P.~G. {Lejeune-Dirichlet}.
\newblock {Vorlesungen \"uber Zahlentheorie. Hrsg. und mit Zus\"atzen versehen
  von R. Dedekind. 4. umgearb. u. verm. Aufl.}
\newblock {Braunschweig. F. Vieweg u. Sohn. XVII + 657 S. \(8^\circ\)}, 1894.

\bibitem{Lie84}
S.~{Lie}.
\newblock {Allgemeine Untersuchungen \"uber Differentialgleichungen, die eine
  continuirliche endliche Gruppe gestatten.}
\newblock {\em {Math. Ann.}}, 25:71--151, 1884.

\bibitem{Macaulay1903}
F.~S. {Macaulay}.
\newblock {Some formulae in eliminations.}
\newblock {\em {Proc. Lond. Math. Soc.}}, 35:3--27, 1903.

\bibitem{Macaulay13}
F.~S. {Macaulay}.
\newblock {On the resolution of a given modular system into primary systems
  including some properties of \textit{Hilbert} numbers.}
\newblock {\em {Math. Ann.}}, 74:66--121, 1913.

\bibitem{Macaulay16}
F.~S. {Macaulay}.
\newblock {The algebraic theory of modular systems.}
\newblock {Cambridge: University press, XIV u. 112 S. $8^{\circ}$.}, 1916.

\bibitem{Macaulay27}
F.~S. {Macaulay}.
\newblock {Some properties of enumeration in the theory of modular systems.}
\newblock {\em {Proc. Lond. Math. Soc. (2)}}, 26:531--555, 1927.

\bibitem{Mansfield96}
E.~L. Mansfield.
\newblock A simple criterion for involutivity.
\newblock {\em J. London Math. Soc. (2)}, 54(2):323--345, 1996.

\bibitem{Mazliak13}
L.~Mazliak.
\newblock {\em Le Carnet de voyage de Maurice Janet à Göttingen}.
\newblock Collection <<\,Essais\,>>. {\'E}ditions Matériologiques, 01 2013.

\bibitem{MerayRiquier1889}
Ch. {M\'eray} and Ch. {Riquier}.
\newblock {Sur la convergence des d\'eveloppements des int\'egrales ordinaires
  d'un syst\`eme d'\'equations diff\'erentielles totales.}
\newblock {\em {Ann. Sci. \'Ec. Norm. Sup\'er. (3)}}, 6:355--378, 1889.

\bibitem{MerayRiquier1890}
Ch. {M\'eray} and Ch. {Riquier}.
\newblock {Sur la convergence des d\'eveloppements des int\'egrales ordinaires
  d'un syst\`eme d'\'equations diff\'erentielles partielles.}
\newblock {\em {Ann. Sci. \'Ec. Norm. Sup\'er. (3)}}, 7:23--88, 1890.

\bibitem{Noether1921}
E.~{Noether}.
\newblock Idealtheorie in ringbereichen.
\newblock {\em Mathematische Annalen}, 83:24--66, 1921.

\bibitem{Noether23}
E.~{Noether}.
\newblock Eliminationstheorie und allgemeine {I}dealtheorie.
\newblock {\em Math. Ann.}, 90(3-4):229--261, 1923.

\bibitem{Pfaff15}
J.~F. {Pfaff}.
\newblock {Allgemeine Methode, partielle Differentialgleichungen zu integrieren
  (1815). Aus dem Lateinischen \"ubersetzt und herausgegeben von {\it G.
  Kowalewski.}.}
\newblock {84 S. $8^{\text {vo}}$ (Ostwalds Klassiker No. 129)}, 1902.

\bibitem{Picard1916}
\'E. {Picard}.
\newblock {L'histoire des sciences et les pr\'etentions de la science
  allemande.}
\newblock {Paris: Perrin, 49 S. \(16^\circ\)}, 1916.

\bibitem{Riquier93}
Ch. Riquier.
\newblock {De l'existence des int\'egrales dans un syst\`eme diff\'erentiel
  quelconque.}
\newblock {\em {Ann. Sci. \'Ec. Norm. Sup\'er. (3)}}, 10:65--86, 123--150,
  167--181, 1893.

\bibitem{Riquier97}
Ch. Riquier.
\newblock Sur les syst\`emes diff\'erentiels les plus g\'en\'eraux.
\newblock {\em Ann. Sci. \'Ecole Norm. Sup. (3)}, 14:99--108, 1897.

\bibitem{Riquier10}
Ch. {Riquier}.
\newblock {Les syst\`emes d'\'equations aux d\'eriv\'ees partielles.}
\newblock {XXVII - 590 p. Paris, Gauthier-Villars.}, 1910.

\bibitem{Riquier28}
Ch. {Riquier}.
\newblock {La m\'ethode des fonctions majorantes et les syst\`emes
  d'\'equations aux d\'eriv\'ees partielles.}
\newblock {Paris: Gauthier-Villars (M\'emorial des sciences math\'ematiques,
  fasc. 32). 63 p. }, 1928.

\bibitem{Schwarz92}
F.~Schwarz.
\newblock An algorithm for determining the size of symmetry groups.
\newblock {\em Computing}, 49(2):95--115, 1992.

\bibitem{Seiler10}
W.~M. Seiler.
\newblock {\em Involution}, volume~24 of {\em Algorithms and Computation in
  Mathematics}.
\newblock Springer-Verlag, Berlin, 2010.
\newblock The formal theory of differential equations and its applications in
  computer algebra.

\bibitem{Serret1849}
J.~A. {Serret}.
\newblock {Cours d'alg\`ebre sup\'erieure.}
\newblock {Paris, Bachelier, Imprimeur-Libraire}, 1849.

\bibitem{Serret1928}
J.-A. {Serret}.
\newblock {Cours d'alg\`ebre sup\'erieure. I, II. 7 \'ed.}
\newblock {Paris: Gauthier-Villars. XIII, 647, XII, 694 p.}, 1928.

\bibitem{Shirshov62}
A.~I. Shirshov.
\newblock Some algorithmic problems for lie algebras.
\newblock {\em Sib. Mat. Zh.}, 3:292--296, 1962.

\bibitem{Thomas37}
J.-M. {Thomas}.
\newblock {Differential systems.}
\newblock {IX + 118 p. New York, American Mathematical Society (American
  Mathematical Society Colloquium Publications Vol. XXI)}, 1937.

\bibitem{Tresse94}
A.~{Tresse}.
\newblock {Sur les invariants diff\'erentie1s des groupes continus de
  transformations.}
\newblock {\em {Acta Math.}}, 18:1--88, 1894.

\bibitem{van-der-Waerden1930}
B.~L. {van der Waerden}.
\newblock {\em {Moderne Algebra. Bd. I. Unter Benutzung von Vorlesungen von
  \textit{E. Artin} und \textit{E. Noether}.}}, volume~23.
\newblock Springer, Berlin, 1930.

\bibitem{Weber00}
E.~{Weber}.
\newblock {Vorlesungen \"uber das Pfaff'sche Problem und die Theorie der
  partiellen Differentialgleichungen erster Ordnung.}
\newblock {Leipzig: B. G. Teubner. XI + 622 S. gr. $8^\circ$}, 1900.

\end{thebibliography}
\end{small}

\vfill

\begin{footnotesize}
\auteur{Kenji Iohara}{iohara@math.univ-lyon1.fr}
{Universit\'{e} de Lyon,
Universit\'{e} Lyon 1, \\
CNRS, Institut Camille Jordan UMR 5208,\\
F-69622 Villeurbanne, France}

\bigskip

\auteur{Philippe Malbos}{malbos@math.univ-lyon1.fr}
{Universit\'{e} de Lyon,
Universit\'{e} Lyon 1, \\
CNRS, Institut Camille Jordan UMR 5208,\\
F-69622 Villeurbanne, France}
\end{footnotesize}

\vspace{3cm}

\begin{small}---\;\;\today\;\;-\;\;\hhmm\;\;---\end{small} \hfill
\end{document}